\documentclass[twocolumn]{autart}

\usepackage{url}
\usepackage{amssymb}
\usepackage{graphicx}
\usepackage{algorithm,algorithmic}	
\usepackage{tikz}
\usepackage{mathtools}
\usepackage[compress]{cite}
\usepackage{amsmath}
\usepackage{dsfont}

\newtheorem{problem}{Problem}
\newtheorem{proposition}{Proposition}

\newtheorem{remark}{Remark}

\newtheorem{assumption}{Assumption}
\newtheorem{theorem}{Theorem}
\newtheorem{corollary}{Corollary}
\newtheorem{example}{Example}

\newcommand{\diag}{\mathop{\text{diag}}}
\newcommand{\argmin}{\mathop{\text{argmin}}}
\newcommand{\argmax}{\mathop{\text{argmax}}}

\renewcommand{\vec}{\mathrm{vec}}

\newcommand{\trace}{\mathrm{Tr}}

\newcommand{\supp}{\mathrm{supp}}

\renewcommand{\footnoterule}{
  \kern -2pt
  \hrule width 0.3\textwidth height .5pt
  \kern 2pt
}

\newcommand{\IEEEQED}{~\rule[-1pt]{5pt}{5pt}\par\medskip}

\newenvironment{IEEEproof}{{\bf Proof:\ }}{ \hfill \IEEEQED}

\begin{document}

\begin{frontmatter}

\title{Ensuring Privacy with Constrained Additive Noise by Minimizing Fisher Information\vspace{-.2in}}

\thanks[footnoteinfo]{An early version of this paper has appeared at the 56th IEEE Conference on Decision and Control, 2017~\cite{CDCpaper_2017}. The work of F.~Farokhi was supported by the McKenzie Fellowship from the University of Melbourne, a competitive grant (MyIP: ID6874) from Defence Science and Technology Group (DSTG), and the VESKI Victoria Fellowship from the Government of Victoria. The work of H. Sandberg was supported by the EU CHIST-ERA project COPES and the Swedish Civil Contingencies Agency through the CERCES project.}

\author[Farhad1,Farhad2]{Farhad~Farokhi}\ead{ffarokhi@unimelb.edu.au},\ead{farhad.farokhi@data61.csiro.au} 
\author[Kalle]{Henrik Sandberg}\ead{hsan@kth.se}

\address[Farhad1]{CSIRO's Data61, Australia}
\address[Farhad2]{Department of Electrical and Electronic Engineering at the University of Melbourne, Australia}
\address[Kalle]{Department of Automatic Control, KTH Royal Institute of Technology, Sweden\vspace{-.2in}}


\begin{abstract}
The problem of preserving the privacy of individual entries of a database when responding to linear or nonlinear queries with constrained additive noise is considered. For privacy protection, the response to the query is systematically corrupted with an additive random noise whose support is a subset or equal to a pre-defined constraint set. A measure of privacy using the inverse of the trace of the Fisher information matrix is developed. The Cram\'{e}r-Rao bound relates the variance of any estimator of the database entries to the introduced privacy measure.
The probability density that minimizes the trace of the Fisher information (as a proxy for maximizing the measure of privacy) is computed. An extension to dynamic problems is also presented. Finally, the results are compared to the differential privacy methodology.\vspace{-.2in}
\end{abstract}

\end{frontmatter}

\section{Introduction}
%\subsection{Motivation}
The constant state of connectedness has enabled the use of new technologies, such as participatory sensing and big-data analysis. These technologies can vastly improve the efficiency of existing infrastructures with little investment. This has come at the price of the erosion of privacy in the society. An example of this lack of privacy is the potential use of data from smart electrical meters by adversaries, such as criminals, advertising agencies, and governments, for monitoring the presence and the activities of occupants~\cite{mcdaniel2009security}. Other examples can include the use of detailed travel data for traffic estimation in intelligent transportation systems~\cite{hoh2006enhancing}, privacy violations caused by sharing information in distributed control systems~\cite{huang2014cost}, and privacy concerns in cloud computing and control~\cite{pearson2010privacy}. These concerns have therefore motivated an urgent need for creating appropriate mechanisms that can protect the privacy of the individuals whose information is stored in various databases. Specifically, it is of interest to provide responses to queries from policy makers on the aggregated data as accurately as possible while not leaking the private information of the individuals.

%\subsection{Contributions}
To combat these problems, the problem of preserving the privacy of individual entries of a database using a constrained additive noise is considered in this paper. It is assumed that anyone, including an adversary, can submit queries to a (trusted) server possessing the entire database. The server returns a response to the query that is systematically corrupted by an additive noise whose support is a subset, or equal, to a desired constraint set. The Cram\'{e}r-Rao bound~\cite[p.\,169]{cramerraotheorem} is then used to relate the variance of the estimation error of unbiased estimators of the database by adversaries  from the provided responses to the trace of the inverse of the Fisher information matrix. 

We start with finding a maximizer of the inverse of the trace of the Fisher information matrix. This optimization problem is nonconvex. It is proved that finding the probability density function of the noise boils down to solving a nonlinear partial differential equation, which is complex in general, even for the simplest case (of linear queries and when the probability density of the noise is independent of the content of the database). Thus, we opt for maximizing a lower bound of the inverse of the trace of the Fisher information matrix, which is the inverse of the trace of the Fisher information matrix scaled by a constant. This is equivalent to minimizing the trace of the Fisher information matrix. This optimization problem is proved to be convex. It is shown that the optimal noise distribution can be calculated by solving a linear partial differential equation (that can be sometimes further simplified with the aid of separation of variables). These results are subsequently generalized to the case where the support set of the noise distribution is unbounded. Noting that, for unbounded sets, the solution to the problem is to add a noise with infinite variance (because the trace of the Fisher information matrix can be pushed to zero), the need for minimizing the trace of the Fisher information (for ensuring privacy) is balanced with the quality of the response (captured by the variance of the additive noise). It can be shown that the Gaussian noise is optimal if the noise is not constrained. These results are demonstrated on three illustrative examples involving smart meter privacy, and computing the average and variance of private databases.

The problem formulation and parts of the results are extended to dynamic estimation problems, where the initial condition of the system is assumed to be the variable that needs to be kept private. This is motivated by a traffic estimation problem in which the initial condition of the system (modelling the vehicle) corresponds to the location of driver's house, which is private.

Finally, the optimal privacy-preserving policies of this paper are compared with differentially-private policies (specifically, the Laplace mechanism) and the optimal privacy-preserving policies when using mutual information as a measure of privacy. It is observed that the optimal policies in the unconstrained-noise formulation are also $(\epsilon,\delta)$-differentially private. Further, the optimal policies in this paper coincide with the optimal privacy-preserving policies when using mutual information as a measure of privacy for the unconstrained case, thus inheriting strong information-theoretic guarantees.

%\subsection{Related Studies}
A common approach for ensuring that the privacy of participants in large databases (or rather the content of the entries of the database owned by those participants) is the application of differential privacy~\cite{dwork2008differential,Dwork2011, le2014differentially,han2014differentially, huang2014cost,7431982}. Those studies often advocate the addition of noises with slow-decaying probability density functions to the response of the queries submitted to the server. This is done so that an adversary cannot accurately infer the private information of the individuals stored in the database. The Laplace noise is frequently utilized in the differential privacy literature; see, e.g.,~\cite{Dwork2014AFD26930522693053,Dwork2011}. 
However, other noise distributions are also common for achieving the differential privacy, or variants thereof. This has prompted various studies to seek the optimal noise distribution for differential privacy~\cite{soria2013optimal,geng2014optimal,7039713}.

Although several information-theoretic interpretations of differential privacy have been presented~\cite{Alvim2012,Mir2013,makhdoumi2013privacy}, the available literature does not offer an operational meaning for the concept as well as a systematic approach for setting the differential privacy parameter (except a broad sweep). Further, in practice, it may not be possible to use an additive noise with infinite support as the noise might need to satisfy certain constraints, e.g., it must belong to a bounded set for smart metering~\cite{farokhisandberg2016}. 

Application of differential privacy in control systems has also gained attention recently. Differential private filtering is discussed in~\cite{le2014differentially}, where releasing filtered signals while respecting privacy of  user data streams is considered. Distributed control of multi-agent systems in the presence of privacy constraints is studied in~\cite{huang2014cost}. An attainable lower bound on the entropy of output is presented in the case where an additive noise is used to ensure differential privacy for discrete-time systems~\cite{wang2014entropy}.  Differential privacy in the specific case of consensus-seeking algorithms is also considered in~\cite{huang2012differentially, nozari2015differentially,katewa2015protecting, mo2016privacy,duan2015privacy,le2015privacy}. A thorough review of these results can be found in a recent tutorial paper~\cite{7798915}. 

Differential privacy has also been successfully utilized in numerical optimization. In~\cite{7431982}, parameters of individual constraints in resource allocation problems is kept private. Preserving the privacy of decisions and cost functions  in distributed optimizations is investigated in~\cite{hale2015differentially, nozari2016differentially}. Applications of differential privacy can also be found in other related problems, such as machine learning~\cite{friedman2010data, rubinstein2012learning}, mechanism design~\cite{mcsherry2007mechanism, nissim2012approximately}, and transportation systems~\cite{dong2015differential, kargl2013differential}.

Recently, several studies have used mutual information (or entropy) and the least mean square estimation error as measures of privacy~\cite{5622047,tan2013increasing,yao2013privacy, tanaka2015sdp,akyol2015privacy,farokhiNecsys2016,farokhi2015quadratic}. Similar to differential privacy literature, most results based on mutual information do not provide an intuitive or interpretable bound on the statistics of the estimation error by the adversary (with the exception of~\cite{TanakaSandberg2017} which uses rate distortion theory to get an interpretable bound on the performance of the adversary). They also require \textit{a priori} assumptions on the distribution of the database, which might not be available in practice due to complexity and scale of the database. The privacy results using the least mean square estimation error also restrict the behaviour of the adversary and assume the underlying random variables are Gaussian, which might not be the case in practice as well. 

In~\cite{farokhisandberg2016}, Fisher information is utilized as a measure of privacy to design privacy-preserving charging policies for batteries in households with smart meters. In this paper, those results are extended to develop a general framework using Fisher information as a measure of privacy. This paper extends~\cite{farokhisandberg2016} in the following ways.
The results of~\cite{farokhisandberg2016}, in the language of this paper, are about releasing all the entries of the database, i.e., the query submitted to the server is an identity function. That very special and restrictive case does not even cover linear queries, let alone providing the optimal privacy-preserving policy for non-linear queries and dynamic estimation.

Finally, it is worth mentioning that the statistics community has previously used Fisher information as a measure of privacy~\cite{anderson1977efficiency,nayak2009unified}. However, in those studies, minimizing the Fisher information to obtain privacy-preserving policies over the set of density functions whose support sets is appropriately constrained  is not discussed.

%\subsection{Paper Outline}
The rest of the paper is organized as follows. The problem formulation and some preliminary results are presented in Section~\ref{sec:problem}. The optimal privacy-preserving probability density functions for the additive constrained noise are developed in Section~\ref{sec:all_results}. These results are then generalized to dynamic estimation problems in Section~\ref{sec:extension}. The relationships between the presented framework and the existing results in the literature are discussed in more depth in Section~\ref{sec:discussions}.  Finally, Section~\ref{sec:conc} concludes the paper and presents viable avenues for future work.

\section{Background and Problem Formulation} \label{sec:problem}
Let $x\in\mathcal{X}\subseteq\mathbb{R}^n$ be a variable that should be kept private, i.e., a database that is possessed by a (trusted) server. This data is only available to the server. In what follows, the vector $x$ is assumed to be deterministic and fixed, i.e., no prior is required nor available. Anyone, including an adversary, can submit queries of the form $f(x)$ with $f:\mathbb{R}^{n}\rightarrow \mathbb{R}^{m}$ to the server. The server in return provides a noisy response to the query given by
\begin{align} \label{eqn:linear_query}
y=f(x)+w,
\end{align}
where $w\in\mathcal{W}(x)\subseteq\mathbb{R}^m$ denotes the noise that the server adds to the data to protect the entries of the private variable $x$. The following standing assumption is made.

\begin{assumption}
$f$ is continuously differentiable. 
\end{assumption}

In this paper, two specific families of constraints on the support set of the additive noise $\mathcal{W}(x)$ are considered:
\begin{itemize}
\item[(\textit{i})] the set $\mathcal{W}(x)$ is independent of $x$; or
\item[(\textit{ii})] the set $\mathcal{W}(x)$ takes the special form of $\{-f(x)\}\oplus\mathcal{Y}$ for some set $\mathcal{Y}\subseteq\mathbb{R}^m$.
\end{itemize}
Here, $\mathcal{C}\oplus\mathcal{B}$ denotes the set $\{a+b\,|\,a\in\mathcal{C},b\in\mathcal{B}\}$ for any two sets $\mathcal{C}$ and $\mathcal{B}$. The family of constraints following the form of (\textit{i}) models the case where the noise itself is constrained, e.g., the additive noise should be positive or bounded. However, the family of constraints in (\textit{ii}) captures the case where the output $y$ must be constrained inside the set $\mathcal{Y}$. In what follows, we use $\mathcal{W}$ to denote either of these families of constraints on the support set of the additive noise.

Note that the server has the intention to respond as accurately as possible to the query $f(x)$ as this typically corresponds to statistical properties (e.g., mean) of the database, which are valuable to, e.g., policy makers. However, it does not want the entries of $x$ (the private data of the people) to be released online nor estimated. Examples of applications where the server wants to provide accurate answers to the submitted queries while keeping the entries of the database hidden can be found in~\cite{Zhou2011,erlich2014routes, dankar2013practicing}. 

The server's policy (which is the object of interest in this paper) is the probability density function $\gamma(\cdot|x):\mathcal{W}\rightarrow \mathbb{R}_{\geq 0}$ of the noise $w$, where $\mathbb{R}_{\geq 0}:=\{x\in\mathbb{R}\,|\,x\geq 0\}$. This implies that
$\mathbb{P}\{w\in\mathcal{W}'\,|\,x\}=\int_{w'\in\mathcal{W}'}\gamma(w'|x)\mathrm{d}w'$ for any Lebesgue-measurable set $\mathcal{W}'\subseteq \mathcal{W}$. In this paper, it is desired to seek a policy $\gamma$ that makes the problem of inferring the private variable~$x$ difficult (according to an appropriate measure described below). The set of all admissible policies $\Gamma$ is restricted according to the following standing assumption.

\begin{assumption} \label{assum:1} (\textit{i})~$\gamma(w|x)$ is such that $\mathbb{P}\{w\in\mathbb{R}^m\setminus\mathcal{W}\,|\,x\}=0$, (\textit{ii}) $\gamma(w|x)$ is twice continuously differentiable in $(w,x)$ over $\mathcal{W}\times\mathcal{X}$, and (\textit{iii}) $\gamma(w|x)=0$ for all $w\in\partial \mathcal{W}$.
\end{assumption}

Assumption~\ref{assum:1}~(\textit{i}) ensures that, with probability one, the noise is restricted to the set $\mathcal{W}$, i.e., $\mathbb{P}\{w\in\mathcal{W}\,|\,x\}=1$. This is to ensure that the constraints on the noise or the output are satisfied almost surely. Assumption~\ref{assum:1}~(\textit{ii}) is required for the use of the Cram\'{e}r-Rao bound as well as the use of results from calculus of variations for finding the optimal probability density function.  Finally, as observed later in the paper, Assumption~\ref{assum:1}~(\textit{iii}) is necessary for the  Cram\'{e}r-Rao bound (see Proposition~\ref{tho:1}). The latter part of this assumption is satisfied if the set $\mathcal{W}$ is unbounded (as, in the limit, a probability density function is always zero otherwise it does not integrate to one). Note that, for bounded constraint sets $\mathcal{W}$, the set of probability density functions that are zero over $\partial \mathcal{W}$ can approximate any probability density function arbitrarily closely. This is proved in~\cite{CDCpaper_2017}.

Under the aforementioned policy of the server, the probability density of $y$ for a given $x$ is then equal to
\begin{align}
p(y|x)&=\gamma(y-f(x)|x),\;\forall y\in \mathcal{Y}=\{f(x)\}\oplus \mathcal{W}.
\end{align}
Further, for any continuously differentiable function $g:\mathbb{R}^n\rightarrow\mathbb{R}$, the notation $\partial g(x)/\partial x$ is used to denote a column vector containing the partial derivatives of the function. For any multivariate function $g(x)$, $G(x)$ denotes its Jacobian, i.e., a matrix with the element in the $i$-th row and the $j$-th column being equal to $\partial g_i(x)/\partial x_j$. Before stating the next preliminary result, the  Fisher information matrix $\mathcal{I}(x)\in\mathbb{R}^{n\times n}$ is defined as
\begin{align}
\mathcal{I}(x)
&\hspace{-.03in}=\hspace{-.03in}\int_{y\in \mathcal{Y}}\hspace{-.1in}
p(y|x)\bigg[\frac{\partial \log(p(y|x))}{\partial x} \bigg]\bigg[\frac{\partial \log(p(y|x))}{\partial x} \bigg]^\top \hspace{-.05in}\mathrm{d}y.
\end{align}
The following results immediately follows from the use of the Cram\'{e}r-Rao bound providing a lower bound on the adversary's estimation error of the private variable independent of the policy. 

\begin{proposition} \label{tho:1} Under Assumption~\ref{assum:1}, for any unbiased estimate of $x$ denoted by $\hat{x}(y)$, it holds that
\begin{align}
\mathbb{E}\{\|x-\hat{x}(y)\|_2^2\}
&\geq \trace(\mathcal{I}(x)^{-1}),
\end{align}
if $\mathcal{I}(x)$ is invertible. Furthermore, for any biased estimate of $x$ denoted by $\hat{x}(y)$ such that $\mathbb{E}\{\hat{x}(y)\}=g(x)$, it holds that
\begin{align}
\mathbb{E}\{\|x-\hat{x}(y)\|_2^2\}
\geq& \trace(G(x)^\top\mathcal{I}(x)^{-1}G(x))\nonumber\\&+\|x-g(x)\|_2^2,
\end{align}
if $\mathcal{I}(x)$ is invertible.
\end{proposition}

\begin{IEEEproof} The proof 
follows from the use of the Cram\'{e}r-Rao bound~\cite[p.\,169]{cramerraotheorem}. 
\end{IEEEproof}

Here, it is desirable to find a policy $\gamma$ for the server that makes estimation of $x$ as difficult as possible. This can be pursued through multiple avenues. When an unbiased estimator of the database exists, following Proposition~\ref{tho:1}, in order to make the task of inferring about $x$ from the measurement $y$ difficult, the trace of the inverse of the Fisher information matrix should be maximized. However, if $m<n$, there may not exist any unbiased estimator of the database $x$ because there are more unknowns $n$ than measurements $m$. In this case, the goal becomes to maximize $\trace(G(x)^\top\mathcal{I}(x)^{-1}G(x))$. It can be shown that
\begin{align}
\trace(G(x)^\top\mathcal{I}(x)^{-1}G(x))
& =
\trace(\mathcal{I}(x)^{-1}G(x)G(x)^\top) \nonumber \\
&\hspace{-.03in}\geq\hspace{-.03in}
\trace(\mathcal{I}(x)^{\hspace{-.02in}-1})\lambda_{\min}(G(x)G(x)^{\hspace{-.03in}\top}\hspace{-.03in}),\label{eqn:lowerboundeig}
\end{align}
where, for any matrix, $\lambda_{\min}(\cdot)$ denotes its smallest eigenvalue. This inequality shows that, even if a biased estimator is utilized, the trace of the inverse of the Fisher information matrix can be maximized to preserve the privacy of the entries of the database, albeit if $G(x)$ assumes full row rank for all $x\in\mathcal{X}$; otherwise the lower bound in~\eqref{eqn:lowerboundeig} is zero and maximizing $\trace(\mathcal{I}(x)^{-1})$ does not result in any tangible privacy guarantee. Therefore, the following assumption is made.

\begin{assumption} \label{assum:rank} $G(x)$ is full row rank for all $x\in\mathcal{X}$.
\end{assumption}

The sensibility of this problem formulation relies on the validity of Assumption~\ref{assum:rank}, which is in general difficult, if not impossible, to check (see Example~\ref{example:average} for a case in which this assumption can be easily checked). Assumption~\ref{assum:rank} is not strictly-speaking necessary, at least for the convexified problem discussed later (see Problem~\ref{problem:general:2}) as the problem formulation can be alternatively motivated by using a worst-case privacy guarantee. Consider the case where, for any $1\leq i\leq n$, the server aims at protecting the content of $x_i$ even if all other entries of the database $x_{-i}=(x_1,\dots,x_{i-1},x_{i+1},\dots,x_{n})$ are leaked (i.e., this is a worst-case analysis for privacy protection). This is to ensure that even if the owners of all the other entries of the database are colluding with the adversary, they cannot extract the private data of an individual (to a reasonable extent). To do so, let $\hat{x}_i(y,x_{-i})$ denote an unbiased estimator of $x_i$ based on $y$ for a given $x_{-i}$. By fixing $x_{-i}$ as knowns, Proposition~\ref{tho:1} can be used to deduce that
\begin{align}
\mathbb{E}\{\|x_i-\hat{x}_i(y,x_{-i})\|_2^2\}
&\geq 1/\mathcal{I}_i(x),
\end{align}
where 
\begin{align*}
\mathcal{I}_i(x)
&=\int_{y\in \mathcal{Y}}
p(y|x)\bigg[\frac{\partial \log(p(y|x))}{\partial x_i} \bigg]\bigg[\frac{\partial \log(p(y|x))}{\partial x_i} \bigg]^\top \mathrm{d}y
\end{align*}
with $p(y|x)=\gamma(w-f(x_i,x_{-i})|x_i,x_{-i})$. Hence,
\begin{align*}
\min_{i} \mathbb{E}\{\|x_i-\hat{x}_i(y,x_{-i})\|_2^2\}
&\geq 
\min_{i} (1/\mathcal{I}_i(x))\\
&= 1/(\max_{i} \mathcal{I}_i(x)).
\end{align*}
Further, note that
\begin{align*}
\max_{i} \mathcal{I}_i(x)
&\leq \sum_{i=1}^n \mathcal{I}_i(x)\\
&=\int_{y\in \mathcal{Y}}
p(y|x)\sum_{i=1}^n\bigg[\frac{\partial \log(p(y|x))}{\partial x_i} \bigg]\\
&\hspace{.5in}\times\bigg[\frac{\partial \log(p(y|x))}{\partial x_i} \bigg]^\top \mathrm{d}y\\
%&=\int_{y\in \mathcal{Y}}
%p(y|x)\trace\bigg(\bigg[\frac{\partial \log(p(y|x))}{\partial x} \bigg]\\
%&\hspace{.5in}\times\bigg[\frac{\partial \log(p(y|x))}{\partial x} \bigg]^\top\bigg) \mathrm{d}y\\
&=\trace\bigg(\int_{y\in \mathcal{Y}}
p(y|x)\bigg[\frac{\partial \log(p(y|x))}{\partial x} \bigg]\\
&\hspace{.5in}\times\bigg[\frac{\partial \log(p(y|x))}{\partial x} \bigg]^\top \mathrm{d}y\bigg)\\
&=\trace(\mathcal{I}(x)).
\end{align*}
Hence, it can be deduced that
\begin{align}
\min_{i} \mathbb{E}\{\|x_i-\hat{x}_i(y,x_{-i})\|_2^2\}
&\geq 1/\trace(\mathcal{I}(x)).
\end{align}
This shows that, even for this stronger notion of privacy (in the absence of Assumption~\ref{assum:rank}), minimizing $\trace(\mathcal{I}(x))$ (instead of maximizing $\trace(\mathcal{I}(x)^{-1})$) provides a reasonable solution.

Before presenting the problem formulation, it should be noted that maximizing $\trace(\mathcal{I}(x)^{-1})$ is not well-defined when $x$ is not known \textit{a priori}. Therefore, instead, we maximize the cost function
\begin{align} \label{eqn:mathcalJ}
\overline{\mathcal{J}}:=\int_{x\in\mathcal{X}} \trace(\mathcal{I}(x)^{-1})p(x)\mathrm{d}x,
\end{align}
where $p:\mathcal{X}\rightarrow\mathbb{R}_{\geq 0}$ is a weight associated with $x$ such that $p(x)\neq 0$ for some set with non-zero Lebesgue measure (since otherwise $\overline{\mathcal{J}}$ is identical to zero). Note that $p(x)$ is not a prior but a weight that captures how difficult the server wants the estimation of the database to become for a given $x$. Note that it can always be assumed that $\int_{x\in\mathcal{X}}p(x)\mathrm{d}x=1$. This is without loss of generality as $p(x)$ can be always scaled by $\int_{x\in\mathcal{X}}p(x)\mathrm{d}x>0$ to achieve the equality. This is of course a design parameter. 

\begin{problem} \label{problem:general:1} Find  $\gamma^*\in\argmax_{\gamma\in\Gamma}\overline{\mathcal{J}}$.
\end{problem}

\begin{remark}[Well-Defined Problem Formulation]
Problem~\ref{problem:general:1} is well-defined if the support set of the additive noise $\mathcal{W}$ is bounded. If this is not the case, the optimal solution is to push the variance of the additive noise to infinity (as that pushes $\trace(\mathcal{I}^{-1})$ to infinity). Note that if $\mathcal{W}$ is bounded, $\trace(\mathcal{I}^{-1})$ is also bounded. Hence, for the case where $\mathcal{W}$ is unbounded, an additional term capturing the quality of the provided response by the server can be included in the utility function to ensure the existence of non-trivial and implementable solutions. This case is investigated later in the paper.
\end{remark}

\begin{remark}[Existence of Solutions] Investigating existence of solutions to Problem~\ref{problem:general:1} is a daunting task due to the nature of the set $\Gamma$. In the remainder of this section, a convex approximation of this problem is presented. In this case, the necessary condition for optimality is also sufficient. This allows us to formulate the problem of finding (sub)optimal privacy-preserving policies as solving a linear partial differential equation with Dirichlet boundary conditions. For the relaxed problem, therefore, the existence of an optimal solution can be cast as the existence of solutions to a linear partial differential equation. For some cases, it is possible to find a solution to the partial differential equation satisfying all the boundary conditions (thus guaranteeing the existence of solutions constructively). 
\end{remark}

\begin{remark}[Data Processing Inequality] Mutual information is often lauded as a measure of privacy due to data processing inequality, i.e., additional manipulation of the transmitted messages based on ones private information can only decrease the amount of the leaked information. This is also considered a  beneficial property of the differential privacy, that is, additional manipulations of the outcomes of a differentially private policy cannot decreases the privacy guarantees of the process~\cite{dwork2014algorithmic}. This property also holds for the Fisher information~\cite{zamir1998proof} pointing to that $\trace(\mathcal{I}(x)^{-1})$ can only be increased upon further manipulations of the transmitted response $y$.
\end{remark}

\begin{remark}[Side Channel Information] Privacy studies, including the presented framework based on Fisher information as well as those based on mutual information and differential privacy, are often fragile to admitting side channel information, e.g., measurements of the private variable already available. In this paper's problem formulation, the primary goal is to maximize $\mathbb{E}\{\|x-\hat{x}(y)\|_2^2\}$ (or its lower bound given by the Fisher information) which clearly states that the adversary's estimator $\hat{x}(y)$ is only a function of $y$. The setup however can be generalized following the same line of reasoning to maximize $\mathbb{E}\{\|x-\hat{x}(y,z)\|_2^2\}$, where $z$ models the side channel information. If the nature of the side-channel information is not fully know, the server can attempt at maximizing $\min_{\upsilon(z|x)}\mathbb{E}\{\|x-\hat{x}(y,z)\|_2^2\}$, where $\upsilon(z|x)$ denoting the conditional probability of the side-channel information given the state can vary over restricted set of density functions. This is an interesting avenue for future research.
\end{remark}

Define the relaxed cost function
\begin{align} \label{eqn:mathcalJ_relaxed}
\mathcal{J}:=\int_{x\in\mathcal{X}} \trace(\mathcal{I}(x))p(x)\mathrm{d}x.
\end{align}
Since maximizing $\overline{\mathcal{J}}$ is difficult in general, the lower bound $n^2\mathcal{J}^{-1}$ can be maximized to achieve a sub-optimal solution; the inequality follows from the application of Proposition~\ref{prop:lowerbound} in Appendix~\ref{proof:prop:lowerbound}. Before stating the  problem formulation, it is beneficial to note that $\mathcal{J}$ is in fact a convex function of~$\gamma$ (over a subset of $\Gamma$) and, thus, the new problem formulation is more computationally feasible. This proof follows from the convexity of the trace of the Fisher information matrix. Note that the proof of the convexity of the Fisher information for scalar random variables is widely available; see for example~\cite[pp.\,80-81]{huber1981} and \cite{Elke1995}. The proof for the multivariate case can also be found in~\cite{farokhisandberg2016}. Define the support of a density function $\gamma(\cdot|x)$ as $\supp(\gamma(\cdot|x)):=\{w\in\mathcal{W}\,|\,\gamma(w|x)>0\}$. 

\begin{assumption} \label{assum:4} $\gamma(w|x)$ is such that $\mathcal{W}\setminus \supp(\gamma(\cdot|x))$ has a zero Lebesgue measure for all $x\in\supp(p)$.
\end{assumption}

Let $\overline{\Gamma}\subseteq\Gamma$ be the set of conditional probability density functions satisfying Assumptions~\ref{assum:1} and~\ref{assum:4}. The following proposition shows that $\mathcal{J}$ is a convex function over $\overline{\Gamma}$ (which is a convex set) and thus  stationarity conditions are necessary and sufficient for finding a minimizer of $\mathcal{J}$ over $\overline{\Gamma}$. 

\begin{proposition} \label{prop:convex:measurement} $\mathcal{J}$ is a convex function of the density function $\gamma$ over the set of density functions $\overline{\Gamma}$.
\end{proposition}

\begin{IEEEproof} The proof is similar to the one in~\cite{farokhisandberg2016} and is thus omitted for the sake of space.
\end{IEEEproof}

Similarly, Proposition~\ref{prop:convex:measurement} motivates us to search for the minimizer of $\mathcal{J}$ over the set of all density functions $\gamma(\cdot|x)$ that are at most over a measure-zero set equal to zero in $\mathrm{int}(\mathcal{W})$, which satisfies the definition of $\overline{\Gamma}$.

\begin{problem} \label{problem:general:2} Find $\gamma^*\in\argmin_{\gamma\in \overline{\Gamma}} \mathcal{J}$.
\end{problem}

As stated earlier, Problems~\ref{problem:general:1} and~\ref{problem:general:2} are both only well-defined for bounded noise support sets $\mathcal{W}$. For unbounded sets $\mathcal{W}$, the solution to both problems is to add a noise with infinite variance. By selecting, for instance, a Gaussian noise with an arbitrarily large variance, the trace of the Fisher information matrix can be pushed towards zero (and the trace of its inverse towards infinity by Proposition~\ref{prop:lowerbound} in Appendix~\ref{proof:prop:lowerbound}). To overcome this problem, the quality of the provided response by the server needs to be balanced with the guaranteed privacy. To do so, a measure of quality of the response can be defined as follows:
\begin{align*}
\mathcal{Q}
&=\int_{x\in\mathcal{X}}p(x)\mathbb{E}\{\|y-f(x)\|_2^2\,|\,x\}\mathrm{d}x\\
&=\int_{x\in\mathcal{X}}\int_{w\in\mathcal{W}}w^\top w \gamma(w|x)p(x)\mathrm{d}x.
\end{align*}
The smaller $\mathcal{Q}$ is, the better the quality of the provided response to the query is. Note that the measure of quality~$\mathcal{Q}$ is a convex function of $\gamma$ since it is linear in the probability density function. Now, the problem formulation can be revised for unbounded constraint set $\mathcal{W}$. 

\begin{problem} \label{problem:not_meas:3} Find $\gamma^*\in\argmin_{\gamma\in\overline{\Gamma}} \mathcal{J}+\varrho \mathcal{Q}$,  where $\varrho>0$ is a constant balancing the need for preserving privacy with the quality of the provided response by the server.
\end{problem}

An alternative problem formulation of the following form can be presented in which a hard constraint on the quality of the response is enforced.

\begin{problem} \label{problem:not_meas:4} Find $\gamma^*\in\argmin_{\gamma\in\overline{\Gamma}:\mathcal{Q}\leq \vartheta}\mathcal{J}$, where $\vartheta>0$ denotes the upper bound on performance degradation caused by the additive noise.
\end{problem}

In light of~\cite{jeyakumar1990zero}, Problems~\ref{problem:not_meas:4} and~\ref{problem:not_meas:3} are equivalent in the sense that, for all $\varrho>0$, there exists $\vartheta>0$ such that the solution of Problem~\ref{problem:not_meas:4} is a solution of Problem~\ref{problem:not_meas:3} and \textit{vice versa}.

With the problem formulations at hand, we are ready to calculate the optimal policy of the server. This is the topic of the next section. 

\section{Privacy-Preserving Policy}
\label{sec:all_results}
In this section, the solutions of the previously-stated problem formulations is presented. We start with Problem~\ref{problem:general:1}. In this case, the solution is given in the following the theorem for linear queries of the form $f(x)=Cx$ with $C\in\mathbb{R}^{m\times n}$ and the case where $\gamma(w|x)$ is independent of $x$, e.g., when measurements of the database are not available. Note that, in this case, the Fisher information matrix can be simplified to
\begin{align*}
\mathcal{I}(x)
&=\int_{y\in \{Cx\}\oplus \mathcal{W}}
\gamma(y-Cx)\bigg[\frac{\partial \log(\gamma(y-Cx))}{\partial x} \bigg]\\
&\hspace{1in}\times\bigg[\frac{\partial \log(\gamma(y-Cx))}{\partial x} \bigg]^\top \mathrm{d}y\\
%&=\int_{y\in \{Cx\}\oplus \mathcal{W}}
%\gamma(y-Cx)C^\top\bigg[\frac{\partial \log(\gamma(w))}{\partial w} \bigg]_{w=y-Cx}\\
%&\hspace{1in}\times \bigg[\frac{\partial \log(\gamma(w))}{\partial w} \bigg]_{w=y-Cx}^\top C \mathrm{d}y\\
&=\int_{w\in \mathcal{W}}
\gamma(w)C^\top\hspace{-.04in}\bigg[\frac{\partial \log(\gamma(w))}{\partial w} \bigg] \hspace{-.04in}\bigg[\frac{\partial \log(\gamma(w))}{\partial w} \bigg]^\top C \mathrm{d}w.
\end{align*}
In this case, $\mathcal{I}(x)$ is no longer a function of~$x$ and is thus simply denoted by $\mathcal{I}$. Further, it should be noted that $\overline{\mathcal{J}}=\trace(\mathcal{I}^{-1})$ and $\mathcal{J}=\trace(\mathcal{I})$. Therefore, the solution of the problem in this case becomes independent of the choice of $p(x)$. Proving this results for more general cases creates several complications without providing more insight. The subsequent results are however proved for nonlinear queries and general policies.  The next theorem presents a necessary condition for the solution of non-convex optimization problem in Problem~\ref{problem:general:1}. 

\begin{theorem} \label{tho:problem:general:1}
Let $\gamma^*(w)$ denote a solution of Problem~\ref{problem:general:1} for linear queries of the form $f(x)=Cx$ over the set of probability density functions that are independent of $x$. Then, it satisfies the following conditions with $u(w)=\sqrt{\gamma^*(w)}$ and some constant $\mu\in\mathbb{R}$:
\begin{align} \label{eqn:tho:problem:general:1}
\begin{cases}
\trace(\mathcal{I}^{-2} C^\top D^2u(w)C)+\mu u(w)=0, & w\in\mathcal{W},\\
u(w)=0, & w\in\partial\mathcal{W},\\
u(w)\neq 0, & w\in\mathrm{int}\mathcal{W},\\
\int_{w\in \mathcal{W}} u(w)^2\mathrm{d}w=1.
\end{cases}
\end{align}
\end{theorem}

\begin{IEEEproof} The proofs are moved to the appendices to avoid interrupting the flow of the presentation. See Appendix~\ref{proof:tho:problem:general:1}.
\end{IEEEproof}

\begin{remark}
In Theorem~\ref{tho:problem:general:1}, $\mu$ denotes the Lagrange multiplier associated with the equality constraint $\int_{w\in \mathcal{W}} u(w)^2\mathrm{d}w=1$ (to  ensure that $\gamma(w)=u(w)^2$ is a probability density function). The conditions in~\eqref{eqn:tho:problem:general:1} are equivalent to the Karush--Kuhn--Tucker (KKT) conditions for the infinite-dimensional optimization problem in Problem~\ref{problem:general:1}. In the rest of the paper, for some specific cases, the value of the multiplier is calculated explicitly. However, in general, the value of the multiplier should be iteratively changed (e.g., using the methods in primal-dual optimization) to find the appropriate value.
\end{remark}

\begin{remark}[Complexity of the Solution]
\label{remark:complex} Note that the partial differential equation in~\eqref{eqn:tho:problem:general:1} is nonlinear because $\mathcal{I}$ in $\trace(\mathcal{I}^{-2} C^\top D^2u(w)C)$ is a function of $u(w)$. Further, Theorem~\ref{tho:problem:general:1} only provides a necessary condition, i.e., the solution of Problem~\ref{problem:general:1} satisfies~\eqref{eqn:tho:problem:general:1} but the reverse does not necessarily hold. As mentioned earlier, these difficulties stem from the complexity of maximizing $\trace(\mathcal{I}^{-1})$, which is a non-concave cost function. 
\end{remark}

Following Remark~\ref{remark:complex}, in the remainder of the paper, the relaxed formulation in Problem~\ref{problem:general:2} and its variants are studied. In what follows, $\mathds{1}_n$ denotes the $n$-dimensional vector of ones. If the dimension $n$ is clear from the context, $\mathds{1}$ is used instead of $\mathds{1}_n$. Further, $F(x)$ denotes the Jacobian of the multivariate function $f(x)$, i.e., a matrix with the element in the $i$-th row and the $j$-th column being equal to $\partial f_i(x)/\partial x_j$. The following theorem provides necessary and sufficient conditions for capturing the solution of Problem~\ref{problem:general:2}.

\begin{theorem} \label{tho:problem:general:2} The solution of Problem~\ref{problem:general:2} is given by $\gamma^*(w|x)=u(w,x)^2$, where $u(w,x)$ satisfies
\begin{align} \label{eqn:tho:nonlinear_query}
\begin{cases}
\trace\hspace{-.03in}\left(\hspace{-.03in}\begin{bmatrix}
F(x) F(x)^\top & F(x) \\ F(x)^\top & I
\end{bmatrix}\hspace{-.03in} D^2 u(w,x)\hspace{-.03in}\right) & \\ \hspace{0.2in}+
L(w,x)\hspace{-.05in}
\begin{bmatrix}
\dfrac{\partial \gamma(w|x)}{\partial w} \\[0.8em] \dfrac{\partial \gamma(w|x)}{\partial x}
\end{bmatrix}
\hspace{-.05in}+\hspace{-.03in}\mu(x) u(w,x)\hspace{-.03in}=\hspace{-.03in}0, &  \hspace{-.1in}w\in\mathcal{W},\\
u(w,x)=0, & \hspace{-.1in}w\in\partial\mathcal{W},\\
u(w,x)\neq 0, & \hspace{-.1in} w\in \mathrm{int}\mathcal{W},\\
\int_{w\in \mathcal{W}} u(w,x)^2\mathrm{d}w=1,
\end{cases}
\end{align}
for some mapping $\mu:\mathcal{X}\rightarrow\mathbb{R}$ and 
\begin{align} \label{eqn:definition_L}
L(w,x)\hspace{-.03in}:=\hspace{-.04in}\begin{bmatrix}
\hspace{-.01in}\dfrac{1}{p(x)}\dfrac{\partial p(x)}{\partial x}^{\hspace{-.03in}\top}\hspace{-.06in}F(x)^\top\hspace{-.04in}+\hspace{-.03in}\mathds{1}^\top\hspace{-.03in} D^2f(x) &\hspace{.1in} \dfrac{1}{p(x)}\dfrac{\partial p(x)}{\partial x}^{\hspace{-.03in}\top} \hspace{-.01in}
\end{bmatrix}\hspace{-.04in}.
\end{align}
Further, all solutions (if multiple) satisfying~\eqref{eqn:tho:nonlinear_query} exhibit the same cost.
\end{theorem}

\begin{IEEEproof} 
See Appendix~\ref{proof:tho:problem:general:2}.
\end{IEEEproof}

\begin{remark} Note that the partial differential equation in~\eqref{eqn:tho:nonlinear_query} is often classified as a semi-linear equation in the sense that it is linear in the partial derivatives (thus it is a linear differential equation) but the coefficients can be potentially non-linear functions of the independent variable (which makes it ``space'' varying). Solving these partial differential equations, in general, is a complex task and out of the scope of this paper. In what follows, the partial differential equation in~\eqref{eqn:tho:nonlinear_query} is solved for all scalar (potentially nonlinear) queries in Corollary~\ref{cor:nonlinear=linear}.
\end{remark}

In what follows, the partial differential equation in~\eqref{eqn:tho:nonlinear_query} is solved for three special cases explicitly to gain some insight into the structure of the solution of Problem~\ref{problem:general:2}. 

\begin{example}[Smart Meter Privacy]
\label{example:smart}
Let the energy consumption of household $k\in\mathcal{N}:=\{1,\dots,n\}$ in a neighbourhood with $n\in\mathbb{N}$ houses be denoted by $x_k\in\mathbb{R}_{\geq 0}$. These variables can be aggregated into a vector to get $x=[x_1 \; \cdots \; x_n]^\top$. Assume that the data is stored on an online server so that it can be studied by policy makers and academics. To avoid unintentionally leaking the private details of the customers, the trusted server in possession of the data adds an appropriate noise to the outcome of the queries to which it responds. Note that the additive noise should be somewhat restricted as arbitrary corruptions might render the data useless or unrealistic. For instance, adding a large negative noise can mask possibly wasteful behaviour of the participants or, in extreme, can transform their combined consumption negative, which might be physically impossible if they do not generate any power using renewable energies. In the case where the server decides to publicly release a sanitized version of the consumption data, it must be assumed that $f(x)=x$ and $\mathcal{W}=[\underline{w},\overline{w}]^n$ with constants $0<\underline{w}\leq \overline{w}<+\infty$. The optimal privacy preserving policy for this example is provided in the following corollary.\hfill $\diamondsuit$
\end{example}

\begin{corollary} \label{cor:problem2:optimal}
Let $n=m$, $p(x)=p$, $f(x)=x$, and $\mathcal{W}=[\underline{w},\overline{w}]^m$ with $-\infty<\underline{w}\leq \overline{w}<+\infty$. The solution of Problem~\ref{problem:general:2} is given by
\begin{align*}
\gamma^*(w|x)\hspace{-.04in}=\hspace{-.04in}&\left(\hspace{-.02in}\frac{2}{\overline{w}\hspace{-.02in}-\hspace{-.02in}\underline{w}}\hspace{-.02in}\right)^{\hspace{-.04in}m}\hspace{-.02in}\prod_{i=1}^m\hspace{-.02in}\cos^2\hspace{-.04in}\left(\frac{\pi}{\overline{w}\hspace{-.02in}-\hspace{-.02in}\underline{w}}\hspace{-.04in}\left(w_i-\frac{\overline{w}\hspace{-.02in}+\hspace{-.02in}\underline{w}}{2}\right)\hspace{-.04in}\right)\hspace{-.04in}\mathds{1}_{w\in\mathcal{W}}.
\end{align*}
\end{corollary}

\begin{IEEEproof} See Appendix~\ref{proof:cor:problem2:optimal}.
\end{IEEEproof}

It is expected that, by increasing $\overline{w}-\underline{w}$ in Example~\ref{example:smart}, it becomes easier to preserve the privacy of the server because the response can be buried deeper within the noise. Although adding more noise can improve the privacy, it also reduces the quality of the provided response by the server. The trade-off between privacy and quality is captured in the following corollary.

\begin{corollary} \label{cor:2} For the optimal policy in Corollary~\ref{cor:problem2:optimal}, the following statements hold:
\begin{itemize}
\item[(\textit{i})] $\mathbb{E}\{\|x-\hat{x}(y)\|_2^2\}\geq \trace(\mathcal{I}^{-1})=\kappa (\overline{w}-\underline{w})^2$ with a constant $\kappa>0$ for any unbiased estimate of $x$ denoted by $\hat{x}(y)$;
\item[(\textit{ii})] $\mathcal{Q}=m(\overline{w}-\underline{w})^2 (2\pi^2-3)/(6 \pi^2)\approx 0.2827m(\overline{w}-\underline{w})^2$.
\end{itemize}
\end{corollary}

\begin{IEEEproof} See Appendix~\ref{proof:cor:2}.
\end{IEEEproof}

Corollary~\ref{cor:2} shows that
$\mathcal{Q}=(6 \pi^2\kappa/(2\pi^2-3))\trace(\mathcal{I}^{-1})$
for the optimal policy-preserving policy. Changing any parameter, such as $\overline{w}$ or $\underline{w}$, that may increase the privacy guarantees inevitably degrades the quality of the response. Therefore, as expected, privacy and quality are conflicting criteria.

\begin{example}[Computing Weighted Average] 
\label{example:average}
A common query that users want to perform on large databases is to calculate the weighted average of  private variables. However, the server in possession of the data does not want the original data to be extracted from the reported weighted average. This operation can be modelled by $f(x)=Cx$, where $C\in\mathbb{R}^{1\times n}$ is such that $C\mathds{1}_n=1$.  The adversary, for instance, upon using the least square approach to find an estimate of the database, gets
$\hat{x}(y)=C^\dag y,$ 
where $C^\dag$ denotes the Moore--Penrose pseudoinverse of $C$, defined as $C^\top(C C^\top)^{-1}$. Assumption~\ref{assum:rank} holds in this example if $C$ has a full row rank, which is a reasonable assumption as otherwise the measurements are not independent. Recalling that $C$ is not a full column rank matrix (i.e., it is only full row rank and $m<n$), $\hat{x}(y)$ is not an unbiased estimator as
$\mathbb{E}\{\hat{x}(y)\}
=C^\dag \mathbb{E}\{y\}=C^\dag Cx.$
For this estimator, it can be deduced that
\begin{align*}
\mathbb{E}\{\|x-\hat{x}(y)\|_2^2\}
=&\|(I-C^\dag C)x\|_2^2+\mathbb{E}\{\|C^\dag w\|_2^2\}\\
%=&\|(I-C^\dag C)x\|_2^2+\trace(V_{ww}(C C^\top)^{-1})\\
\geq &\|(I-C^\dag C)x\|_2^2+\trace(\mathcal{I}(x)^{-1}(C C^\top)^{-1}),
\end{align*}
where the last inequality follows from the Cram\'{e}r-Rao bound in Proposition~\ref{tho:1}. Finally, it is also worth saying that, in this case, $\trace(\mathcal{I}(x)^{-1}(C C^\top)^{-1})=\trace(\mathcal{I}(x)^{-1})(C C^\top)^{-1}$ because $m=1$ (and thus $(C C^\top)^{-1}\in\mathbb{R}$). 
This problem clearly fits the framework presented in this paper. Here, the set $\mathcal{W}$ can be either bounded or unbounded based on the situation. The optimal privacy preserving policy for this example is provided in the following corollary.\hfill $\diamondsuit$
\end{example}

\begin{corollary} \label{cor:problem2:averaging}
Let $n\in\mathbb{N}$, $m=1$, $p(x)=p$, $f(x)=Cx$, and $\mathcal{W}=[\underline{w},\overline{w}]$ with $-\infty<\underline{w}\leq \overline{w}<+\infty$. If $C\neq 0$, the solution of Problem~\ref{problem:general:2} is given by
\begin{align} \label{eqn:solution:cor:problem2:averaging}
\gamma^*(w|x)=\frac{2}{\overline{w}-\underline{w}}\cos^2\left(\frac{\pi}{\overline{w}-\underline{w}}\left(w-\frac{\overline{w}+\underline{w}}{2}\right)\right)\mathds{1}_{w\in\mathcal{W}}.
\end{align}
\end{corollary}

\begin{IEEEproof} See Appendix~\ref{proof:cor:problem2:averaging}.
\end{IEEEproof}

Example~\ref{example:average} evidently fits the criteria of Corollary~\ref{cor:problem2:averaging}. It is interesting to note that the choice of the weights in $C$ (so long as at least one of them is non-zero) is irrelevant. The next corollary captures the effect of the weight $p(x)$ on the optimal privacy-preserving policy.

\begin{corollary} \label{cor:nonuniformp} Let $n=m=1$, $f(x)=x$, and $\mathcal{W}=[\underline{w},\overline{w}]$ with $-\infty<\underline{w}\leq \overline{w}<+\infty$. The solution of Problem~\ref{problem:general:2} is given by
\begin{align*}
\gamma^*(w|x)=&c(x)\exp\bigg(-\frac{p'(x)}{p(x)}(w+x)\bigg)\\
&\times\cos^2\left(\frac{\pi}{\overline{w}-\underline{w}}\left(w-\frac{\overline{w}+\underline{w}}{2}\right)\right)\mathds{1}_{w\in\mathcal{W}}
\end{align*}
where
\begin{align*}
c(x)=\bigg[\int_{\underline{w}}^{\overline{w}}&
\exp \bigg(-\frac{p'(x)}{p(x)}(w+x)\bigg)\\&\times\cos^2\left(\frac{\pi}{\overline{w}-\underline{w}}\left(w-\frac{\overline{w}+\underline{w}}{2}\right)\right)\mathrm{d}w
\bigg]^{-1}.
\end{align*}
\end{corollary}

\begin{IEEEproof}
See Appendix~\ref{proof:cor:nonuniformp}.
\end{IEEEproof}

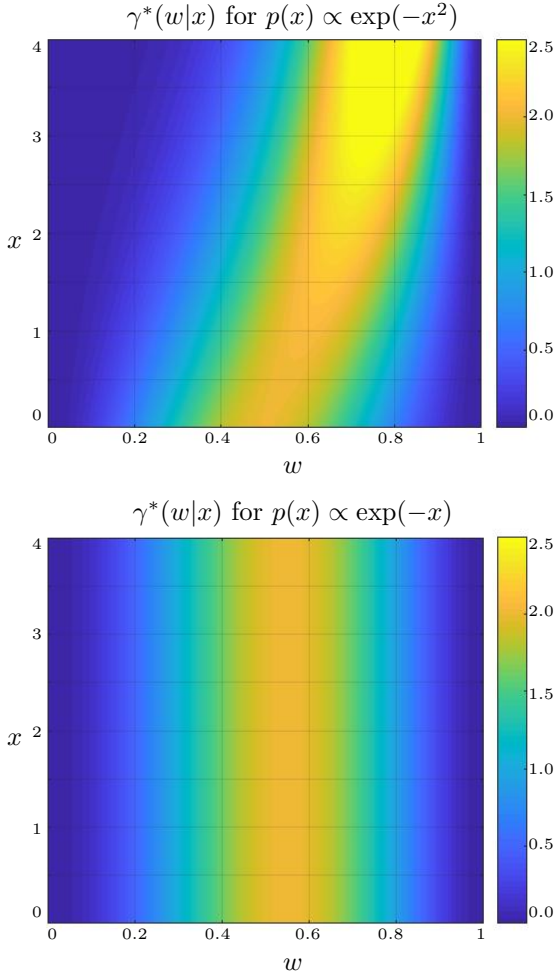
\begin{figure}
\begin{tikzpicture}
\node[] at (0,0) {\includegraphics[width=1\linewidth]{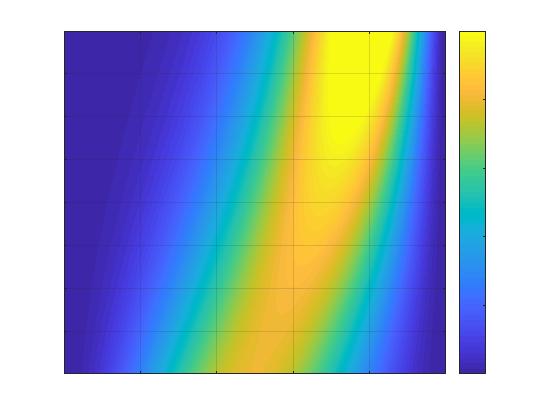}};
\node[] at (-3.2,-2.6) {\tiny $0$};
\node[] at (-2.1,-2.6) {\tiny $0.2$};
\node[] at (-1.0,-2.6) {\tiny $0.4$};
\node[] at (+0.2,-2.6) {\tiny $0.6$};
\node[] at (+1.3,-2.6) {\tiny $0.8$};
\node[] at (+2.5,-2.6) {\tiny $1$};
\node[] at (+0.0,-3.0) {$w$};
\node[] at (-3.4,-2.3) {\tiny $0$};
\node[] at (-3.4,-1.2) {\tiny $1$};
\node[] at (-3.4,+0.1) {\tiny $2$};
\node[] at (-3.4,+1.4) {\tiny $3$};
\node[] at (-3.4,+2.6) {\tiny $4$};
\node[] at (-3.7,-0.0) {$x$};
\node[] at (+3.3,-2.3) {\tiny $0.0$};
\node[] at (+3.3,-1.4) {\tiny $0.5$};
\node[] at (+3.3,-0.4) {\tiny $1.0$};
\node[] at (+3.3,+0.6) {\tiny $1.5$};
\node[] at (+3.3,+1.7) {\tiny $2.0$};
\node[] at (+3.3,+2.6) {\tiny $2.5$};
\node[] at (+0.0,+3.0) {$\gamma^*(w|x)$ for $p(x)\propto\exp(-x^2)$};
\end{tikzpicture}
\begin{tikzpicture}
\node[] at (0,0) {\includegraphics[width=1\linewidth]{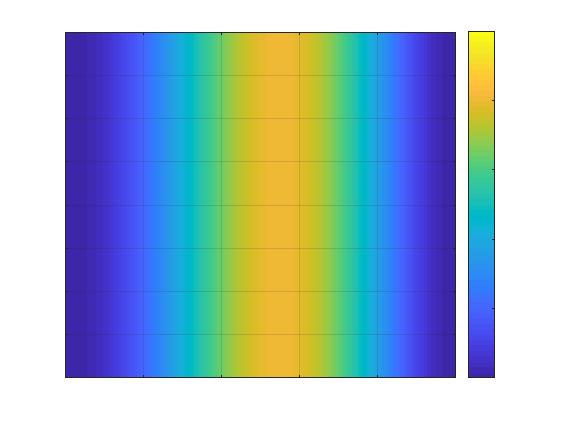}};
\node[] at (-3.2,-2.6) {\tiny $0$};
\node[] at (-2.1,-2.6) {\tiny $0.2$};
\node[] at (-1.0,-2.6) {\tiny $0.4$};
\node[] at (+0.2,-2.6) {\tiny $0.6$};
\node[] at (+1.3,-2.6) {\tiny $0.8$};
\node[] at (+2.5,-2.6) {\tiny $1$};
\node[] at (+0.0,-3.0) {$w$};
\node[] at (-3.4,-2.3) {\tiny $0$};
\node[] at (-3.4,-1.2) {\tiny $1$};
\node[] at (-3.4,+0.1) {\tiny $2$};
\node[] at (-3.4,+1.4) {\tiny $3$};
\node[] at (-3.4,+2.6) {\tiny $4$};
\node[] at (-3.7,-0.0) {$x$};
\node[] at (+3.3,-2.3) {\tiny $0.0$};
\node[] at (+3.3,-1.4) {\tiny $0.5$};
\node[] at (+3.3,-0.4) {\tiny $1.0$};
\node[] at (+3.3,+0.6) {\tiny $1.5$};
\node[] at (+3.3,+1.7) {\tiny $2.0$};
\node[] at (+3.3,+2.6) {\tiny $2.5$};
\node[] at (+0.0,+3.0) {$\gamma^*(w|x)$ for $p(x)\propto\exp(-x)$};
\end{tikzpicture}
\caption{\label{fig:2}Optimal privacy-preserving policy for non-uniform weighting functions $p(x)\propto\exp(-x^2)$ (top) and $p(x)\propto\exp(-x)$ (bottom). }
\end{figure}

Figure~\ref{fig:2} illustrates the optimal privacy-preserving policy $\gamma^*(w|x)$ in Corollary~\ref{cor:nonuniformp} versus $w$ and $x$ for the case where the weighting function is $p(x)\propto\exp(-x^2)$ (top) and $p(x)\propto\exp(-x)$ (bottom) when $\underline{w}=0$ and $\overline{w}=1$. For any mappings $f$ and $g$, it is said that $f(x)\propto g(x)$ if there exists constant $c$ such that $f(x)=cg(x)$. For $p(x)=\exp(-x^2)$, the policy is a function of $x$. However, for $p(x)=\exp(-x)$, the policy is independent of $x$; this can be attributed to that the ratio $p'(x)/p(x)=-1$ is not a function of $x$.

\begin{example}[Computing Variance] 
\label{example:variance}
Another common query, beside the average of a set of private data (see Example~\ref{example:average}), is to calculate its statistical variance. This operation can be done by the nonlinear query
\begin{align*}
f(x)
&=\frac{1}{n-1}\sum_{i=1}^n\bigg(x_i-\frac{1}{n}\sum_{j=1}^n x_j\bigg)^2.
\end{align*}
The optimal privacy preserving policy for this example is provided in the following corollary.\hfill $\diamondsuit$
\end{example}

\begin{corollary} \label{cor:nonlinear=linear}
Let $n\in\mathbb{N}$, $m=1$, $p(x)=p$, and $\mathcal{W}=[\underline{w},\overline{w}]$ with $-\infty<\underline{w}\leq \overline{w}<+\infty$. If $f(x)\neq 0$, the solution of Problem~\ref{problem:general:2} is given by~\eqref{eqn:solution:cor:problem2:averaging}.
\end{corollary}

\begin{IEEEproof} See Appendix~\ref{proof:cor:nonlinear=linear}.
\end{IEEEproof}

Corollary~\ref{cor:nonlinear=linear} proves that the for scalar queries, i.e., when $m=1$, the nonlinearity of the query does not change the distribution of the optimal additive noise. Therefore, the optimal noise distributions for the averaging problem in Example~\ref{example:average} and the variance calculation in Example~\ref{example:variance} are the same.

Now, we are ready to study the case where the constraint set is unbounded. We start by solving Problem~\ref{problem:not_meas:3} in the next theorem.

\begin{theorem} \label{tho:nonlinear_query:unbounded} The solution of Problem~\ref{problem:not_meas:3} is given by $\gamma^*(w|x)=u(w,x)^2$, where $u(w,x)$ satisfies
\begin{align} \label{eqn:tho:nonlinear_query:unbounded}
\begin{cases}
\trace\hspace{-.03in}\left(\hspace{-.03in}\begin{bmatrix}
F(x) F(x)^\top & F(x)^\top \\ F(x) & I
\end{bmatrix}\hspace{-.03in} D^2 u(w,x)\hspace{-.03in}\right) & \\
\hspace{.32in}+
L(w,x)\hspace{-.05in}
\begin{bmatrix}
\dfrac{\partial \gamma(w|x)}{\partial w} \\[0.8em] \dfrac{\partial \gamma(w|x)}{\partial x}
\end{bmatrix}
\hspace{-.03in} & \\ \hspace{.32in}+(\mu(x)\hspace{-.03in}-\hspace{-.03in}(\varrho/4)w^\top w) u(w,x)=0, &  \hspace{-.1in}w\in\mathcal{W},\\
u(w,x)=0, & \hspace{-.1in}w\in\partial\mathcal{W},\\
u(w,x)\neq 0, & \hspace{-.1in}w\in \mathrm{int}\mathcal{W},\\
\int_{w\in \mathcal{W}} u(w,x)^2\mathrm{d}w=1,
\end{cases}
\end{align}
for some mapping $\mu:\mathcal{X}\rightarrow\mathbb{R}$ and $L(w,x)$ is defined in~\eqref{eqn:definition_L}. Further, all solutions (if multiple) satisfying~\eqref{eqn:tho:nonlinear_query:unbounded} exhibit the same cost.
\end{theorem}

\begin{IEEEproof} See Appendix~\ref{proof:tho:4}.
\end{IEEEproof}

An explicit solution of Problems~\ref{problem:not_meas:3} and~\ref{problem:not_meas:4} for the case where $\mathcal{W}=\mathbb{R}^m$ is presented in the following corollary. This case is of special interest as the optimal noise can be compared with the density of the noise suggested in the differential privacy literature, i.e., Laplace noise.

\begin{corollary} \label{cor:problem3:averaging}
Let $\mathcal{W}=\mathbb{R}^m$ and $p(x)=p$. For linear queries of the form $f(x)=Cx$ with  full row rank matrix $C$, the solutions of Problems~\ref{problem:not_meas:3} and~\ref{problem:not_meas:4} is given by
\begin{align*}
\gamma^*(w|x)=\frac{1}{\sqrt{(2\pi)^m\det(\Sigma)}}\exp\bigg(-\frac{1}{2}w^\top \Sigma^{-1} w\bigg),
\end{align*}
where $\Sigma=2(CC^\top)^{1/2}/\sqrt{\rho}$ for Problem~\ref{problem:not_meas:3} and $\Sigma=\vartheta(CC^\top)^{1/2}/\trace((CC^\top)^{1/2})$ Problem~\ref{problem:not_meas:4}. 
\end{corollary}

\begin{IEEEproof} See Appendix~\ref{proof:cor:problem3:averaging}.
\end{IEEEproof}

Corollary~\ref{cor:problem3:averaging} simply states that, in this case, the optimal noise is a (multivariate) Gaussian random variable with covariance $(CC^\top)^{1/2}/\sqrt{\rho}$. Note that the optimality of the Gaussian noise when minimizing the Fisher information over noises with unbounded support set is not in itself new~\cite{5753094}; however, its application in privacy-preserving policies has not been considered previously. Clearly, as $\varrho$ increases, i.e., the emphasis on the quality of the response to the query increases, the variance of the noise decreases. In this framework, the optimal noise distribution differs from the Laplace distribution, which is a standard choice in the differential privacy literature~\cite{Dwork2014AFD26930522693053}. 
Evidently, the weighted averaging setup in Example~\ref{example:average} satisfies the conditions for results of Corollaries~\ref{cor:problem3:averaging}. 
For that example, the optimal noise can be simplified into a Gaussian random variable with variance $\vartheta$ (due to the scalar nature of the responses). Note that, for this example, the optimal privacy-preserving policy is again independent of the choice of the matrix $C$.

\section{Extensions to Dynamic Estimation}
\label{sec:extension}
In this section, the problem formulation is extended by considering dynamic estimation problems for linear time-invariant systems. Assume that
\begin{align*}
x[k+1]&=Ax[k],\quad x[0]=x_0,
\end{align*}
where $x[k]\in\mathbb{R}^n$ is the state. The dynamics, for instance, can capture the case where some entries of a database are getting updated in real time or that the database contains the states of a physical system evolving through time (e.g., position and velocity of a vehicle). Assume that the initial state is deterministic and unknown to the adversary. In each time step $k\in\{0,\dots,T\}$ with $T$ denoting the time horizon, the user can submit a query of the form $Cx[k]$ to the server. The server responds by returning 
\begin{align*}
y[k]=Cx[k]+w[k],
\end{align*}
where $w[k]\in\mathbb{R}^m$ is an additive noise introduced by the server to keep the state private. Here, only the case where the probability density function of $w[k]$ is independent of the state (in the past and the future) is considered.  
Also consider the case where the support of the noise density is unrestricted. Evidently, the results can be extended to the case where the support set of the noise is constrained following the same line of reasoning as in Section~\ref{sec:all_results}. 

Define
$w_k:=[
w[0]^\top \;
w[1]^\top \;
\cdots \;
w[k]^\top
]^\top.$
The policy of the server is the probability density function of $w_T$ denoted by $\gamma:\mathbb{R}^{n(T+1)}\rightarrow\mathbb{R}_{\geq 0}$. This is the most general policy that the server can employ. To keep the entries of the database private, the server wants to increase the covariance of the estimation error of the initial condition $\mathbb{E}\{\|\hat{x}_0(y_T)-x_0\|_2^2\}$, where $\hat{x}_0(y_T)$ denotes any unbiased smoothing estimate of the state based on all the received response aggregated into a single vector of the form
$y_T=[
y[0]^\top \;
y[1]^\top \;
\cdots \;
y[T]^\top
]^\top.$ 
Note that
$y_T=\Psi_Tx_0+w_T,$
where $\Psi_T:=[
C^\top \;
(CA)^\top \;
\cdots \;
(CA^T)^\top
]^\top.$
The conditional density of $y_T$ for any $x_0$ is given by
$p(y_T|x_0)=\gamma(y_T-\Psi_Tx_0).$
Now, the Cram\'{e}r-Rao bound (see Proposition~\ref{tho:1}) can used to show that
\begin{align*}
\mathbb{E}\{\|\hat{x}_0(y_T)-x_0\|_2^2\}
&\geq \trace(\mathcal{I}^{-1})\geq (T+1)^2n^2/\trace(\mathcal{I}),
\end{align*}
where
\begin{align*}
\mathcal{I}
&=\int
\hspace{-0.05in}p(y_T|x_0)\bigg[\frac{\partial \log(p(y_T|x_0))}{\partial x_0} \bigg]\hspace{-0.05in}\bigg[\frac{\partial \log(p(y_T|x_0))}{\partial x_0} \bigg]^{\top}\hspace{-0.05in}\mathrm{d}y_T\\
&=\int\hspace{-0.05in}
\frac{1}{\gamma(w_T)}\Psi_T^\top\bigg[\frac{\partial \gamma(w_T)}{\partial w_T} \bigg]\bigg[\frac{\partial \gamma(w_T)}{\partial w_T} \bigg]^{\top}\Psi_T\mathrm{d}w_T
\end{align*}
and the second inequality follows from Proposition~\ref{prop:lowerbound} in Appendix~\ref{proof:prop:lowerbound}.
The quality of the response can also be measured using
\begin{align*}
\mathcal{Q}
&=\mathbb{E}\{\|y_T-\Psi_T x_0\|_2^2\,|\,x_0\}
%&=\mathbb{E}\{\|w_T\|_2^2\}\\
=\int w_T^\top w_T\gamma(w_T)\mathrm{d}w_T.
\end{align*}
The next theorem provides the optimal policy in the sense of Problem~\ref{problem:not_meas:3} for this case.

\begin{theorem}  \label{tho:dynamic_estimation}
Let $\Psi_T^\top \Psi_T$ be invertible. The solution of~\eqref{problem:not_meas:3} for the dynamic estimation setup is given by $\gamma^*$ such that $\mathbb{P}\{w_T=\Psi_T z\}=1$ where $z$ is distributed according to the probability density function
\begin{align*}
p(z)=\frac{1}{\sqrt{(2\pi)^n\det(\Sigma)}}\exp\bigg(-\frac{1}{2}z^\top \Sigma^{-1} z\bigg),
\end{align*}
where $\Sigma=2(\Psi_T^\top \Psi_T)^{-1/2}/\sqrt{\varrho}$.	
\end{theorem}

\begin{IEEEproof} See Appendix~\ref{proof:tho:dynamic_estimation}.
\end{IEEEproof}

Notice that the noise distribution across time is not independently and identically distributed (i.i.d.) across time. The noise at time step $k$ takes the form of $CA^kz$. This means that the server realizes a random variable $z$ based on the probability distribution in Theorem~\ref{tho:dynamic_estimation}. Then it adds $z$ to the initial condition of the system and propagates its effect in all future time steps.

\begin{remark} \label{prop:invertible} The matrix $\Psi_T^\top \Psi_T$, which is in fact equal to the observability Gramian over $\{0,\dots,T\}$, is invertible if the pair $(A,C)$ is observable and $T\geq n$. 
This follows from observability of linear time-invariant systems~\cite[p.\,271]{sontag2013mathematical} and the Cayley-Hamilton theorem~\cite[p.\,141]{1964methods}. 
\end{remark}

\begin{example}[Traffic Crowd-Sensing with Privacy] Consider the case where a vehicle is sharing its position with a remote monitoring station in pursuit of estimating the state of the traffic on a road. The house of the vehicle's owner is on this road (which is conveniently modelled by the real line) at $s_0\in\mathbb{R}$. The vehicle travels on the road with the constant velocity $v_0\in\mathbb{R}$ after leaving the house. The dynamics of the vehicle is given by
\begin{align*}
\begin{bmatrix}
s[k+1] \\
v[k+1] \\
\end{bmatrix}
&=
\begin{bmatrix}
1 & 1 \\ 0 & 1
\end{bmatrix}
\begin{bmatrix}
s[k] \\
v[k] \\
\end{bmatrix}
,\quad 
\begin{bmatrix}
s[0] \\
v[0] \\
\end{bmatrix}
=\begin{bmatrix}
s_0 \\
v_0 \\
\end{bmatrix},
\end{align*}
where $s[k]$ and $v[k]$, respectively, denote the position and the speed of vehicle over time. The user (i.e., the remote monitoring station) is interested in knowing the position of the vehicle; therefore, it submits a query of the form $Cx[k]$ with $C=[ 1 \; 0 ]$. According to Theorem~\ref{tho:dynamic_estimation}, the server (or the vehicle) needs to provide the response $y[k]=Cx[k]+CA^kz,$ where, if $T>2$, $z$ is a zero mean Gaussian random variable with covariance 
\begin{align*}
\Sigma
&=\frac{2}{\sqrt{\rho}}(\Psi_T^\top \Psi_T)^{-1/2}\\
&=\frac{2}{\sqrt{\rho}}
\begin{bmatrix}
T+1 & T(T+1)/2 \\
T(T+1)/2 & T(T+1)(2T+1)/6
\end{bmatrix}^{-1/2}.
\end{align*}
The quality of the response becomes
\begin{align*}
\mathcal{Q}
&=\frac{1}{\sqrt{3\rho}}\bigg(
\sqrt{
(T + 1)
(2T^2 +T+ 6-
\sqrt{\Delta})} \\
&\quad\quad\quad\quad+ 
\sqrt{(T + 1)
( 2T^2+T + 6 + 
\sqrt{\Delta} )}\bigg),
\end{align*}
where $\Delta=4T^4 + 4T^3 + 13T^2 - 12T + 36$. This implies that\footnote{We say $f(x)=\mathcal{O}(g(x))$ if $\lim_{x\rightarrow\infty} |f(x)/g(x)|=c<\infty$.} $\mathcal{Q}=\mathcal{O}(T\sqrt{T}/\sqrt{\rho})$. On the other hand, 
\begin{align*}
\mathbb{E}\{\|\hat{x}_0(y_T)\hspace{-.03in}-\hspace{-.03in}x_0\|_2^2\}
&\hspace{-.03in}=\hspace{-.03in}\frac{4\sqrt{3}}{\sqrt{\rho}}
\bigg(\hspace{-.03in}\frac{1}{\sqrt{
(T + 1)(2T^2+T + 6 - \sqrt{\Delta} )}} \\
&\quad+\frac{1}{\sqrt{(T + 1)(2T^2+T +6+ \sqrt{\Delta})}}\bigg).
\end{align*}
Thus, $\mathbb{E}\{\|\hat{x}_0(y_T)-x_0\|_2^2\}=\mathcal{O}(1/(T\sqrt{T}\sqrt{\rho}))$. This shows the privacy guarantee decreases with $T$. Therefore, to keep the privacy guarantee constant, smaller $\rho$ should be used for larger horizons $T$ (i.e., a lower emphasis on preserving the quality of response must be placed). In fact, $\rho$ should be selected such that it scales according to $1/(T\sqrt{T})$ with $T$. Doing so, the quality of the response $\mathcal{Q}$ depreciates according to $\mathcal{O}(T^3)$. \hfill $\diamondsuit$
\end{example}

\section{Discussion} \label{sec:discussions}
The choice of the Fisher information as a measure of privacy in this paper, motivated by the Cram\'{e}r-Rao bound, ensures that the privacy guarantees are applicable to a wide range of adversaries in contrast to, e.g.,~\cite{farokhi2015quadratic} assuming a least mean square error estimator as the adversary. Further, the Cram\'{e}r-Rao bound provides a clear operational meaning for the measure of privacy. This could be potentially lacking in differential privacy literature~\cite{Dwork2011} and studies using mutual information and entropy as a measure of privacy, e.g.,~\cite{akyol2015privacy}. However, in the presence of a prior distribution for the data and possible correlations between the entries of the database, privacy-preserving methods that do not use this additional information, such as the proposed method in this paper and algorithms relying on differential privacy, can underperform or break down~\cite{kifer2011no}, which is not the case for method relying on mutual information. 

The use of constrained additive noise in this paper sets it apart from other studies in the literature that use an additive noise whose distribution has an infinite support, such as Laplace or Gaussian~\cite{Dwork2011,Dwork2014AFD26930522693053, han2014differentially,huang2014cost, le2014differentially}. In the studies where the optimal noise is investigated, the support of the distribution is most often unrestricted, which again gives rise to the Laplace or Gaussian distributions being the optimal choices~\cite{akyol2015privacy, farokhi2015quadratic,farokhiNecsys2016,7039713}.

In the unconstrained case, the optimal noise distribution minimizing the Fisher information subject to a constraint on the degradation of the quality of the response is proved to be Gaussian in fact\footnote{Note that by changing the measure of the quality of the response (e.g., expection of the norm-1 of the additive noise), one can get other noise density functions.}. This fact is at odds with the differential privacy literature~\cite{Dwork2011}. Therefore, the provided framework with the guarantees provided by use of the Cram\'{e}r-Rao bound is weaker than the differential privacy (both in requirements and guarantees). The use of the Gaussian noise is however known to satisfy a weaker variant of the differential privacy, referred to as $(\epsilon,\delta)$ differential privacy~\cite{le2014differentially, sandberg2015differentially}. Further, noting that the privacy-preserving policy in this paper and that of~\cite{akyol2015privacy} coincide in the unconstrained case, it is easy to see that the presented framework also minimizes the mutual information while having the advantage of providing a better operational meaning for the measure of privacy. This observation can be intuitively explained by the intimate relationship between Fisher information and mutual information~\cite{clarke1990information, rissanen1996fisher}.

This section is finished by exploring the relationship between differential privacy and the proposed optimal noise in more depth for the weighted averaging setup in Example~\ref{example:average} with $\mathcal{X}=[\underline{x},\overline{x}]^n$ for some $-\infty<\underline{x}\leq \overline{x}<+\infty$. The corrupted response of the server is $\epsilon$-differentially private if
$\mathbb{P}\{Cx+w\in\overline{\mathcal{Y}}\}\leq e^\epsilon \mathbb{P}\{Cx'+w\in\overline{\mathcal{Y}}\}$ for all $x,x'\in\mathcal{X}$ that only differ in one element and all Lebesgue-measurable sets $\overline{\mathcal{Y}}$. The following proposition proves that differential privacy can be achieved by an additive Laplace noise. 

\begin{proposition} \label{prop:differential_privacy} For Example~\ref{example:average} with $\mathcal{X}=[\underline{x},\overline{x}]^n$ for some $-\infty<\underline{x}\leq \overline{x}<+\infty$, the corrupted response of the server is $\epsilon$-differentially private if $\gamma(w)=1/(2b)\exp\left(-|w|/b\right)$ with $b=\epsilon/[(\overline{x}-\underline{x})\max_i |c_i|]$, where $c_i$ is the $i$-th entry of $C$.
\end{proposition}

\begin{IEEEproof} See Appendix~\ref{proof:prop:differential_privacy}.
\end{IEEEproof}

For the differentially private noise in Proposition~\ref{prop:differential_privacy}, it can be shown that $\mathcal{Q}=2b^2$. Therefore, under the constraint $\mathcal{Q}=\vartheta$, one can achieve $[(\overline{x}-\underline{x})\sqrt{\vartheta/2}\max_i c_i]$-differential privacy. Further, it can be shown that for differentially private noise in Proposition~\ref{prop:differential_privacy},
\begin{align*}
\mathcal{I}
&=2\int_{0}^\infty 
\frac{CC^\top}{\gamma(w)}\bigg[\frac{\partial \gamma(w)}{\partial w} \bigg]^2 \mathrm{d}w
=\frac{CC^\top}{b^2}=\frac{2CC^\top}{\vartheta}.
\end{align*}
Therefore, for any unbiased estimator of $x$ denoted by $\hat{x}(y)$ under the noise density function in Proposition~\ref{prop:differential_privacy}, $\mathbb{E}\{\|x-\hat{x}(y)\|_2^2\}\geq \|(I-C^\dag C)x\|_2^2+\vartheta/(2(CC^\top)^2).$
However, for the optimal noise in Corollary~\ref{cor:problem3:averaging}, it can be shown that
$\mathcal{I}=CC^\top/\vartheta.$
Thus, for any unbiased estimator of $x$ denoted by $\hat{x}(y)$ under the noise in Corollary~\ref{cor:problem3:averaging}, $
\mathbb{E}\{\|x-\hat{x}(y)\|_2^2\}\geq \|(I-C^\dag C)x\|_2^2+\vartheta/(CC^\top)^2.$ Note that
\begin{align*}
\sup_{x\in [\underline{x},\overline{x}]^n}&\frac{\|(I-C^\dag C)x\|_2^2+\vartheta/(CC^\top)^2}{\|(I-C^\dag C)x\|_2^2+\vartheta/(2(CC^\top)^2)}
%&=1+\sup_{x\in [\underline{x},\overline{x}]^n}\frac{\vartheta/(2(CC^\top)^2)}{\|(I-C^\dag C)x\|_2^2+\vartheta/(2(CC^\top)^2)}
=1+\kappa>1,
\end{align*}
where
$\kappa=1/({1+2(CC^\top)^{2}\max_{x\in [\underline{x},\overline{x}]^n}\|(I-C^\dag C)x\|_2^2/\vartheta}).$ Thus, for the same bound on the quality of response, the optimal noise distribution in Corollary~\ref{cor:problem3:averaging} provides a privacy guarantee that is $1+\kappa$ times stronger, albeit in the sense of the error covariance $\mathbb{E}\{\|x-\hat{x}(y)\|_2^2\}$. 

For the differentially private noise in Proposition~\ref{prop:differential_privacy}, it can be shown that $H(\gamma)=\log_2(e\sqrt{2\vartheta})$. For the optimal noise in Corollary~\ref{cor:problem3:averaging}, $H(\gamma^*)=\log_2(\sqrt{e2\pi\vartheta})$. Interestingly, $H(\gamma^*)\geq H(\gamma)$, which \textit{intuitively} points to the fact that, for some prior distributions on the entries of the database, the optimal noise in Corollary~\ref{cor:problem3:averaging} can potentially provide stronger information theoretic guarantees as well. This is not particularly surprising considering that the optimal noise for privacy in an information theoretic setting is also proved to be Gaussian~\cite{akyol2015privacy}.

\begin{figure}
\centering
\begin{tikzpicture}
\node[] at (0,0) {\includegraphics[width=1.0\linewidth]{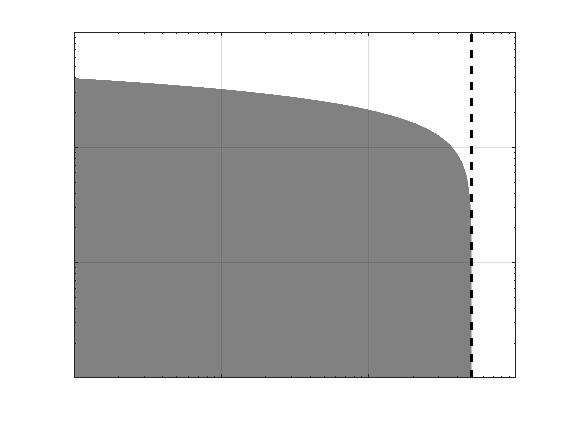}};
\node[rotate=45] at (-3.0,-2.8) {\tiny $10^{-3}$};
\node[rotate=45] at (-0.9,-2.8) {\tiny $10^{-2}$};
\node[rotate=45] at (+1.2,-2.8) {\tiny $10^{-1}$};
\node[rotate=45] at (+3.4,-2.8) {\tiny $10^{0}$};
\node[] at (+0.0,-3.5) {$\delta$};
\node[rotate=45] at (-3.4,-2.3) {\tiny $10^{-3}$};
\node[rotate=45] at (-3.4,-0.7) {\tiny $10^{-2}$};
\node[rotate=45] at (-3.4,+1.0) {\tiny $10^{-1}$};
\node[rotate=45] at (-3.4,+2.6) {\tiny $10^{0}$};
\node[] at (-4.0,-0.0) {$\epsilon$};
\node[] (1) at (+1.0,2.1) {$\delta=1/2$};
\draw[->] (1) -- (+2.7,1.8);
\end{tikzpicture}
\caption{\label{fig:1} Let $\vartheta/((\overline{x}-\underline{x})\max_{i}|c_i|)=1$. The white and gray areas illustrate the regions for which the condition~\eqref{eqn:weak_differential_privacy} is satisfied and is not satisfied, respectively. Note that satisfying~\eqref{eqn:weak_differential_privacy} for $\delta\leq 1/2$ implies that the noise distribution in Corollary~\ref{cor:problem3:averaging}  is $(\epsilon,\delta)$-differentially private.}
\end{figure}
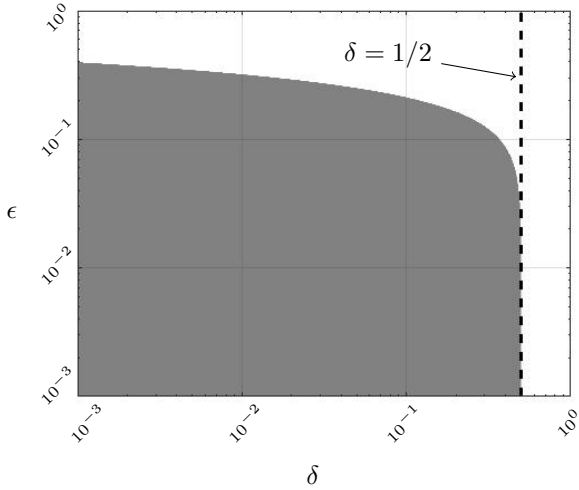

The corrupted response of server is $(\epsilon,\delta)$-differentially private if $\mathbb{P}\{Cx+w\in\overline{\mathcal{Y}}\}\leq e^\epsilon \mathbb{P}\{Cx'+w\in\overline{\mathcal{Y}}\}+\delta$ for all $x,x'\in\mathcal{X}$ that only differ in one element and all Lebesgue-measurable sets $\overline{\mathcal{Y}}$.  The following proposition proves that differential privacy can be achieved by an additive Gaussian noise. 

\begin{proposition} \label{prop:weak_differential_privacy} For Example~\ref{example:average} with $\mathcal{X}=[\underline{x},\overline{x}]^n$ for some $-\infty<\underline{x}\leq \overline{x}<+\infty$, the corrupted response of the server with noise distribution in Corollary~\ref{cor:problem3:averaging}  is $(\epsilon,\delta)$-differentially private if $\delta\leq 1/2$ and 
\begin{align} \label{eqn:weak_differential_privacy}
\vartheta\geq (\overline{x}-\underline{x})\max_{i}|c_i|\bigg(\frac{\sqrt{2\ln(1/(2\delta))}}{\epsilon}+\frac{1}{\sqrt{2\epsilon}} \bigg).
\end{align}
\end{proposition}

\begin{IEEEproof} See Appendix~\ref{proof:prop:weak_differential_privacy}.\end{IEEEproof}

In Figure~\ref{fig:1}, the white and gray areas illustrate the regions for which the condition~\eqref{eqn:weak_differential_privacy} is satisfied and is not satisfied, respectively, when $\vartheta/((\overline{x}-\underline{x})\max_{i}|c_i|)=1$. Note that satisfying~\eqref{eqn:weak_differential_privacy} for $\delta\leq 1/2$ (behind the dashed line) implies that the noise distribution in Corollary~\ref{cor:problem3:averaging} is $(\epsilon,\delta)$-differentially private for the corresponding values of $\epsilon$ and $\delta$.
 
\section{Conclusions and Future Work}\label{sec:conc}
In this paper, the problem of preserving the privacy of individual entries of a database with constrained additive noise is investigated. A measure of privacy using the Fisher information is developed. The optimal probability density function that maximizes the measure of privacy is computed. It is shown that, in some cases, the privacy-preserving policy of the server could potentially be independent of the entries of the database. Further, for scalar queries when the support set of the additive noise is bounded, the nature of the query, e.g., its linearity or content, does not play a role in the optimal additive noise for preserving the privacy of the database. For unconstrained noises, Gaussian distribution seems to be the optimal privacy-preserving policy if the quality of the provided response is measured using the variance of the additive noise employed by the server. Future work can focus on dynamic nonlinear systems.

\bibliographystyle{ieeetr}
\bibliography{ref}

\appendix

\section{A Useful Inequality}
\label{proof:prop:lowerbound}

\begin{proposition} \label{prop:lowerbound} $\overline{\mathcal{J}}\geq n^2\mathcal{J}^{-1}$.
\end{proposition}

\begin{IEEEproof} First, note that
\begin{align*}
\trace(\mathcal{I}^{-1})
=\sum_{i=1}^n \frac{1}{\lambda_i(\mathcal{I})}
\geq \frac{n^2}{\sum_{i=1}^n\lambda_i(\mathcal{I})}=n^2\trace(\mathcal{I})^{-1},
\end{align*}
where the inequality follows from the Jensen's inequality~\cite[p.\,25]{rockafellar2015convex} and the facts that the mapping $z\mapsto 1/z$ is convex over $\mathbb{R}_{\geq 0}$ and $\lambda_i(\mathcal{I})\geq 0$ for all $i$ (since $\mathcal{I}$ is positive semi-definite). Further, it can be shown that
\begin{align*}
\int_{x\in\mathcal{X}} \trace(\mathcal{I}(x)^{-1})p(x)\mathrm{d}x
&\geq n^2\int_{x\in\mathcal{X}} \trace(\mathcal{I}(x))^{-1}p(x)\mathrm{d}x\\
&\geq n^2\bigg(\int_{x\in\mathcal{X}} \trace(\mathcal{I}(x))p(x)\mathrm{d}x\bigg)^{-1},
\end{align*}
where the second inequality follows from the Jensen's inequality and the earlier observation that the mapping $z\mapsto 1/z$ is convex over $\mathbb{R}_{\geq 0}$.
\end{IEEEproof}

\section{Proof of Theorem~\ref{tho:problem:general:1}}
\label{proof:tho:problem:general:1}
The variational derivative of $\trace(\mathcal{I}^{-1})$ needs to be calculated by introducing infinitesimal functional variations in the probability density function $\gamma$. Let $\delta \gamma$ denote the infinitesimal variations in $\gamma$. The variational derivative $(\mathcal{I})_{ij}$, the element of $\mathcal{I}$ in row $i$ and column $j$, is denoted by $\delta (\mathcal{I})_{ij}$. Using Theorem~5.2 in~\cite[p.\,440]{edwards1973advanced}, it can be shown that 
\begin{align*}
\delta (\mathcal{I})_{ij}
=\int_{w\in\mathcal{W}} \bigg[&\hspace{-.03in} \frac{1}{\gamma(w)^2}\bigg(\hspace{-.03in}\sum_{k}c_{ki}\frac{\partial\gamma(w)}{\partial w_k}\bigg)\hspace{-.03in}\bigg(\hspace{-.03in}\sum_{\ell}c_{\ell j}\frac{\partial\gamma(w)}{\partial w_\ell}\bigg)\\&-\frac{2}{\gamma(w)}\sum_{k,\ell}c_{ki}c_{\ell j}\frac{\partial^2\gamma(w)}{\partial w_k\partial w_\ell}\bigg]\delta \gamma(w)\mathrm{d}w,
\end{align*}
where, for any matrix $C$, $c_{ij}$ denotes the entry in the $i$-th row and $j$-th column.
Define $\delta \mathcal{I}$ to be a matrix where the element in row $i$ and column $j$ is equal to $\delta (\mathcal{I})_{ij}$. Using linear algebra, the variational derivative can be rewritten in matrix form as
\begin{align*}
\delta \mathcal{I}
=\int_{w\in\mathcal{W}} \bigg[& \frac{1}{\gamma(w)^2}C^\top\bigg[\frac{\partial\gamma(w)}{\partial w}\bigg]\bigg[\frac{\partial\gamma(w)}{\partial w}\bigg]^\top C\\
&- \frac{2}{\gamma(w)}C^\top D^2\gamma(w)C\bigg]\delta \gamma(w)\mathrm{d}w.
\end{align*}
The variational derivative of $\trace(\mathcal{I}^{-1})$ can be calculated as
\begin{align*}
\delta \trace(\mathcal{I}^{-1})
&=\lim_{\epsilon \rightarrow 0} (\trace((\mathcal{I}+\epsilon \delta \mathcal{I})^{-1})-\trace(\mathcal{I}^{-1}))/\epsilon\\
&=-\trace(\mathcal{I}^{-1}\delta \mathcal{I}\mathcal{I}^{-1})\\
&=-\trace(\mathcal{I}^{-2}\delta \mathcal{I}).
\end{align*}
Now, the Lagrangian can be constructed according to
\begin{align*}
\mathcal{L}
=&\trace(\mathcal{I}^{-1})+\mu \bigg(\int_{w\in \mathcal{W}} \gamma(w)\mathrm{d}w -1\bigg)\\
=&\trace(\mathcal{I}^{-1})+\int_{w\in \mathcal{W}}
\mu \gamma(w) \mathrm{d}b-\mu,
\end{align*}
where $\mu\in\mathbb{R}$ is the Lagrange multiplier corresponding to the equality constraint $\int_{w\in \mathcal{W}} \gamma(w)\mathrm{d}w=1$. The necessary condition for optimality is that the extrema must make the variational derivative of $\mathcal{L}$ equal to zero. As a result,
\begin{align*}
\int_{w\in\mathcal{W}} &\bigg[ \trace\bigg(\mathcal{I}^{-2}\bigg[\frac{1}{\gamma(w)^2}C^\top\bigg[\frac{\partial\gamma(w)}{\partial w}\bigg]\bigg[\frac{\partial\gamma(w)}{\partial w}\bigg]^\top C\\
&-\frac{2}{\gamma(w)}C^\top D^2\gamma(w)C\bigg]\bigg)-\mu \bigg]\delta\gamma(w)\mathrm{d}w=0
\end{align*}
for all $\delta \gamma$. This is only possible if
\begin{align*}
\trace\bigg(\mathcal{I}^{-2}\bigg[\frac{1}{\gamma(w)^2}C^\top&\bigg[\frac{\partial\gamma(w)}{\partial w}\bigg]\bigg[\frac{\partial\gamma(w)}{\partial w}\bigg]^\top C\\
&-\frac{2}{\gamma(w)}C^\top D^2\gamma(w)C\bigg]\bigg)-\mu=0.
\end{align*}
Introducing the change of variable $\gamma(w)=u(w)^2$ results in
\begin{align} \label{eqn:proof:1}
\mu+\frac{4}{u(w)}\trace(\mathcal{I}^{-2}C^\top D^2u(w) C)=0.
\end{align}
If $u(w)\neq 0$ for all $w\in\mathrm{int}(\mathcal{W})$,~\eqref{eqn:proof:1} can be rewritten as
\begin{align*}
\trace(\mathcal{I}^{-2}C^\top D^2u(w) C)+\bar{\mu} u(w)=0,\quad w\in\mathrm{int}(\mathcal{W}),
\end{align*}
where $\bar{\mu}=\mu/4$. However, if $u(w)=0$ for some $w\in\mathrm{int}(\mathcal{W})$, the equality in~\eqref{eqn:proof:1} cannot be satisfied with any $\mu\in\mathbb{R}$.

\section{Proof of Theorem~\ref{tho:problem:general:2}}
\label{proof:tho:problem:general:2}
First, noting that the cost function and the constraint set are convex, the stationary condition (that the variational derivative is equal to zero) is sufficient for optimality.  Further, if multiple density functions satisfy the sufficiency conditions, they all exhibit the same cost. In the rest of the proof, this condition is rewritten in a simpler forms. To do so, note that
\begin{align*}
\mathcal{I}(x)&\hspace{-.03in}=\hspace{-.03in}\int_{y\in \{f(x)\}\oplus \mathcal{W}}
\gamma(y-f(x)|x)\bigg[\frac{\partial \log(\gamma(y-f(x)|x))}{\partial x} \bigg]\\
&\hspace{1in}\times\bigg[\frac{\partial \log(\gamma(y-f(x)|x))}{\partial x} \bigg]^\top \mathrm{d}y\\
%&\hspace{-.03in}=\hspace{-.03in}\int_{y\in \{f(x)\}\oplus \mathcal{W}}
%\gamma(y-f(x)|x)\\
%&\hspace{0.2in}\times\hspace{-.03in}\bigg[F(x)^\top\frac{\partial \log(\gamma(w|x))}{\partial w}\hspace{-.04in}+\hspace{-.04in}\frac{\partial \log(\gamma(w|x))}{\partial x} \bigg]_{w=y-f(x)}\\
%&\hspace{0.2in}\times\hspace{-.03in} \bigg[F(x)^\top\frac{\partial \log(\gamma(w|x))}{\partial w}\hspace*{-.04in}+\hspace{-.04in}\frac{\partial \log(\gamma(w|x))}{\partial x} \bigg]_{w=y-f(x)}^\top \mathrm{d}y\\
&\hspace{-.03in}=\hspace{-.03in}\int_{w\in \mathcal{W}}\hspace{-.2in}
\gamma(w|x)\bigg[F(x)^\top\frac{\partial \log(\gamma(w|x))}{\partial w}+\frac{\partial \log(\gamma(w|x))}{\partial x} \bigg]\\
&\hspace{0.2in}\times \bigg[F(x)^\top\frac{\partial \log(\gamma(w|x))}{\partial w}+\frac{\partial \log(\gamma(w|x))}{\partial x} \bigg]^\top  \mathrm{d}w\\
&\hspace{-.03in}=\hspace{-.03in}\int_{w\in \mathcal{W}}
\frac{1}{\gamma(w|x)}\bigg[F(x)^\top\frac{\partial \gamma(w|x)}{\partial w}+\frac{\partial \gamma(w|x)}{\partial x} \bigg]\\
&\hspace{0.2in}\times \bigg[F(x)^\top\frac{\partial \gamma(w|x)}{\partial w}+\frac{\partial \gamma(w|x)}{\partial x} \bigg]^\top  \mathrm{d}w.
\end{align*}
Thus,
\begin{align*}
\trace(\mathcal{I}(x))
&=\int_{w\in \mathcal{W}}
\frac{1}{\gamma(w|x)}\bigg[F(x)^\top\frac{\partial \gamma(w|x)}{\partial w}+\frac{\partial \gamma(w|x)}{\partial x} \bigg]^\top\\
&\hspace{0.6in}\times \bigg[F(x)^\top\frac{\partial \gamma(w|x)}{\partial w}+\frac{\partial \gamma(w|x)}{\partial x} \bigg] \mathrm{d}w.
\end{align*}
As a result,
\begin{align*}
\mathcal{J}
=&\int_{x\in\mathcal{X}}\int_{w\in \mathcal{W}}
\frac{p(x)}{\gamma(w|x)}\bigg[F(x)^\top\frac{\partial \gamma(w|x)}{\partial w}+\frac{\partial \gamma(w|x)}{\partial x} \bigg]^\top\\
&\hspace{0.6in}\times \bigg[F(x)^\top\frac{\partial \gamma(w|x)}{\partial w}+\frac{\partial \gamma(w|x)}{\partial x} \bigg] \mathrm{d}w\mathrm{d}x.
\end{align*}
Following the result of~\cite{jeyakumar1990zero}, the Lagrangian  can be constructed as
\begin{align*}
\mathcal{L}
=&\int_{x\in\mathcal{X}}\int_{w\in \mathcal{W}}
\frac{p(x)}{\gamma(w|x)}\bigg[F(x)^\top\frac{\partial \gamma(w|x)}{\partial w}+\frac{\partial \gamma(w|x)}{\partial x} \bigg]^\top\\
&\hspace{0.9in}\times \bigg[F(x)^\top\frac{\partial \gamma(w|x)}{\partial w}+\frac{\partial \gamma(w|x)}{\partial x} \bigg] \mathrm{d}w\mathrm{d}x\\
&-\int_{x\in\mathcal{X}}p(x)\mu(x) \bigg(\int_{w\in \mathcal{W}} \gamma(w|x)\mathrm{d}w -1\bigg)\mathrm{d}x\\
=&\int_{x\in\mathcal{X}}\int_{w\in \mathcal{W}}
\hspace{-.05in}p(x)\bigg(\hspace{-.07in}-\hspace{-.03in}\mu(x) \gamma(w|x)\\
&\hspace{0.5in}+
\frac{1}{\gamma(w|x)}\bigg[F(x)^\top\frac{\partial \gamma(w|x)}{\partial w}+\frac{\partial \gamma(w|x)}{\partial x} \bigg]^\top\\
&\hspace{0.9in}\times \bigg[F(x)^\top\frac{\partial \gamma(w|x)}{\partial w}+\frac{\partial \gamma(w|x)}{\partial x} \bigg]\bigg) \mathrm{d}w\mathrm{d}x\\
&\hspace{0.5in}+\int_{x\in\mathcal{X}}\mu(x)p(x)\mathrm{d}x,
\end{align*}
where $\mu:\mathcal{X}\rightarrow\mathbb{R}$ is the Lagrange multiplier corresponding to the equality constraint $\int_{w\in \mathcal{W}} \gamma(w|x)\mathrm{d}w=1$ for all $x\in\supp(p)$. Using Theorem~5.3 in~\cite[p.\,440]{edwards1973advanced}, it can be seen that the extrema must satisfy~\eqref{eqn:extrema:long1} on top of the next page. 
\begin{figure*}
\begin{align}
&\frac{p(x)}{\gamma(w|x)^2}
\begin{bmatrix}
\dfrac{\partial \gamma(w|x)}{\partial w} \\[0.8em] \dfrac{\partial \gamma(w|x)}{\partial x}
\end{bmatrix}^\top\hspace{-.08in}
\begin{bmatrix}
F(x)F(x)^\top & F(x) \\ F(x)^\top & I
\end{bmatrix}\hspace{-.08in}\begin{bmatrix}
\dfrac{\partial \gamma(w|x)}{\partial w} \\[0.8em] \dfrac{\partial \gamma(w|x)}{\partial x}
\end{bmatrix}\hspace{-.06in} +\hspace{-.03in}2\hspace{-.03in}\sum_{i=1}^m \frac{\partial}{\partial w_i}\hspace{-.04in} \left(\hspace{-.04in}\frac{p(x)}{\gamma(w|x)}
e_i^\top\hspace{-.05in}
\begin{bmatrix}
F(x)F(x)^\top & F(x) \\ F(x)^\top & I
\end{bmatrix}\hspace{-.06in} \begin{bmatrix}
\dfrac{\partial \gamma(w|x)}{\partial w} \\[0.8em] \dfrac{\partial \gamma(w|x)}{\partial x}
\end{bmatrix} \hspace{-.04in}\right)\nonumber\\
&+2p(x)\hspace{-.06in}\sum_{i=m+1}^{n+m}\hspace{-.06in} \frac{\partial}{\partial x_i} \left(\frac{1}{\gamma(w|x)}
e_i^\top
\begin{bmatrix}
F(x)F(x)^\top & F(x) \\ F(x)^\top & I
\end{bmatrix} \begin{bmatrix}
\dfrac{\partial \gamma(w|x)}{\partial w} \\[0.8em] \dfrac{\partial \gamma(w|x)}{\partial x}
\end{bmatrix} \right)\nonumber \\
&+2\hspace{-.06in}\sum_{i=m+1}^{n+m}\hspace{-.06in} \frac{\partial p(x)}{\partial x_i} \left(\frac{1}{\gamma(w|x)}
e_i^\top
\begin{bmatrix}
F(x)F(x)^\top & F(x) \\ F(x)^\top & I
\end{bmatrix} \begin{bmatrix}
\dfrac{\partial \gamma(w|x)}{\partial w} \\[0.8em] \dfrac{\partial \gamma(w|x)}{\partial x}
\end{bmatrix} \right)+p(x)\mu(x)
\hspace{-.05in}=\hspace{-.03in}0.
\label{eqn:extrema:long1}
\end{align} 
\hrule
\end{figure*}
Introducing the change of variable $\gamma(w|x)=u(w,x)^2$ results in
\begin{align*}
\mu(x)&+\frac{4}{u(x,w)}\trace\left(\begin{bmatrix}
F(x)F(x)^\top & F(x) \\ F(x)^\top & I
\end{bmatrix} D^2 u(w,x)\right)\\
&+\frac{4}{u(x,w)}\mathds{1}^\top D^2f(x)\frac{\partial u(w,x)}{\partial w}\\
&+\frac{4}{u(x,w)}\frac{1}{p(x)}\frac{\partial p(x)}{\partial x}^\top 
\begin{bmatrix}
F(x)^\top & I
\end{bmatrix}
\begin{bmatrix}
\dfrac{\partial \gamma(w|x)}{\partial w} \\[0.8em] \dfrac{\partial \gamma(w|x)}{\partial x}
\end{bmatrix}
=0,
\end{align*}
for all $w\in\mathrm{int}(\mathcal{W})$ and $x\in\supp(p)$. 
Again, if $u(w,x)\neq 0$ for all $w\in\mathrm{int}(\mathcal{W})$ and $x\in\supp(p)$, it can be deduced that
\begin{align*}
\bar{\mu}(x) u(w,x)&+\trace\bigg(\begin{bmatrix}
F(x)F(x)^\top & F(x) \\ F(x)^\top & I
\end{bmatrix} D^2 u(w,x)\bigg)\\
&+\mathds{1}^\top D^2f(x)\frac{\partial u(w,x)}{\partial w}\\
&+\frac{1}{p(x)}\frac{\partial p(x)}{\partial x}^\top 
\begin{bmatrix}
F(x)^\top & I
\end{bmatrix}
\begin{bmatrix}
\dfrac{\partial \gamma(w|x)}{\partial w} \\[0.8em] \dfrac{\partial \gamma(w|x)}{\partial x}
\end{bmatrix}=0,
\end{align*}
where $\bar{\mu}(x)=\mu(x)/4$. However, if $u(w,x)=0$ for some $w\in\mathrm{int}(\mathcal{W})$ and $x\in\supp(p)$, the equality cannot be satisfied with any $\mu\in\mathbb{R}$. 

\section{Proof of Corollary~\ref{cor:problem2:optimal}}
\label{proof:cor:problem2:optimal}
In this proof, a solution of the form $u(w,x)=u(w)$ and $\mu(x)=\mu$ is sought for the partial differential equation in~\eqref{eqn:tho:nonlinear_query}. 
In this case, the partial differential equation in~\eqref{eqn:tho:nonlinear_query} becomes
\begin{align} \label{eqn:schrodinger}
\begin{cases}
\nabla^2 u(w)+\mu u(w)=0, & w\in\mathcal{W},\\
u(w)=0, & w\in\partial\mathcal{W}.
\end{cases}
\end{align}
This is a special case of the time-independent Schr\"{o}dinger equation. This knowledge can be used to solve the partial differential equation explicitly.
Following~\cite{MohamedAtia2005}, the solution of~\eqref{eqn:schrodinger} is unique. The rest easily follows from showing that the provided density function satisfies the partial differential equation and its boundary conditions. 

\section{Proof of Corollary~\ref{cor:2}}
\label{proof:cor:2}
First, we prove part~(\textit{i}). Allow $u^{(\underline{w},\overline{w})}(w)$ be such that $\gamma(w)=u^{(\underline{w},\overline{w})}(w)^2$ with $\gamma(w)$ denoting the solution of Problem~\ref{problem:general:2} in Corollary~\ref{cor:problem2:optimal} for $\mathcal{W}=[\underline{w},\overline{w}]$. To emphasize the fact that $\mathcal{W}$ is a function of $\underline{w}$ and $\overline{w}$, in this proof, the notation $\mathcal{W}^{(\underline{w},\overline{w})}$ is used.
Therefore,
\begin{align*}
\mathcal{I}
&=\int_{w\in\mathcal{W}^{(\underline{w},\overline{w})}}
\bigg[\frac{\partial u^{(\underline{w},\overline{w})}(w)}{\partial w} \bigg]
\bigg[\frac{\partial u^{(\underline{w},\overline{w})}(w)}{\partial w} \bigg]^\top \mathrm{d}w\\
&=\int_{w\in\mathcal{W}^{(\underline{w},\overline{w})}} \frac{1}{(\overline{w}\hspace{-.04in}-\hspace{-.04in}\underline{w})^{m\hspace{-.02in}+\hspace{-.02in}2}}
\bigg[\frac{\partial u^{(0,1)}(w')}{\partial w'} \bigg]_{w'=(w-\underline{w})/(\overline{w}-\underline{w})} \\&\hspace{1in}\times
\bigg[\frac{\partial u^{(0,1)}(w')}{\partial w'} \bigg]_{w'=(w-\underline{w})/(\overline{w}-\underline{w})}^\top \mathrm{d}w\\
&=\frac{1}{(\overline{w}-\underline{w})^{2}}\underbrace{\hspace{-.04in}\int_{w'\in\mathcal{W}^{(0,1)}} 
\hspace{-.04in}\bigg[\frac{\partial u^{(0,1)}(w')}{\partial w'} \bigg]
\bigg[\frac{\partial u^{(0,1)}(w')}{\partial w'} \bigg]^\top \hspace{-.04in}\mathrm{d}w'}_{:=\kappa}\hspace{-.02in}.
\end{align*}

Now, we prove part~(\textit{ii}). To do so, note that
\begin{align*}
\mathcal{Q}
&=\trace(\mathbb{E}\{ww^\top \})\\
&=\trace(\mathbb{E}\{(w-\mathbb{E}\{w\})(w-\mathbb{E}\{w\})^\top\})+m(\overline{w}+\underline{w})^2/4,
\end{align*}
where the second equality follows from that $\mathbb{E}\{w\}^\top \mathbb{E}\{w\}=m(\overline{w}+\underline{w})^2/4$. Noting that $w_i$ is independent of $w_j$ for the optimal policy in Corollary~\ref{cor:problem2:optimal} results in
\begin{align*}
\mathbb{E}&\{(w-\mathbb{E}\{w\})(w-\mathbb{E}\{w\})^\top\}\\
&=\diag(
\mathbb{E}\{(w_1-\mathbb{E}\{w_1\})^2\},\cdots,\mathbb{E}\{(w_m-\mathbb{E}\{w_m\})^2\}).
\end{align*}
Note that
\begin{align*}
\mathbb{E}\{(w_i\hspace{-.03in}-\hspace{-.03in}\mathbb{E}\{w_i\})^2\}
=&\left(\frac{2}{\overline{w}-\underline{w}}\right)\int_{\underline{w}}^{\overline{w}}\left(w_i-\frac{\overline{w}+\underline{w}}{2}\right)^2\\
&\times\cos^2\left(\hspace{-.03in}\frac{\pi}{\overline{w}-\underline{w}}\hspace{-.03in}\left(w_i-\frac{\overline{w}+\underline{w}}{2}\hspace{-.03in}\right)\hspace{-.03in}\right)\mathrm{d}w_i\\
=&\frac{(\pi^2-6) (\overline{w}-\underline{w})^2}{12 \pi^2}.
\end{align*}

\section{Proof of Corollary~\ref{cor:problem2:averaging}}
\label{proof:cor:problem2:averaging}
In this proof, a solution of the form $u(w,x)=u(w)$ and $\mu(x)=\mu$ is sought for the partial differential equation in~\eqref{eqn:tho:nonlinear_query}. 
In this case, it can be shown that 
$\trace(C^\top D^2u(w) C)
=u''(w)\trace(C^\top C)$.
Therefore, the partial differential equation in~\eqref{eqn:tho:nonlinear_query} becomes the ordinary differential equation $u''(w)+\bar{\mu} u(w)=0$ for $w\in\mathcal{W}$ with the boundary condition that $u(w)=0$ for all $w\in\partial\mathcal{W}$, 
where $\bar{\mu}=\mu/\trace(C^\top C)$. The differential equation $u''(w)+\bar{\mu} u(w)=0$ admits a solution of the form $u(w)=\alpha \cos(\sqrt{\bar{\mu}}(w-\beta))$ where $\alpha,\beta\in\mathbb{R}$ are constants depending on the boundary conditions. It should be ensured that $u(\underline{w})=u(\overline{w})=0$ since $\gamma(w)=u(w)^2$ for $w\in\partial \mathcal{W}$. Thus $\bar{\mu}>0$. Two distinct situations may occur:
\begin{itemize}
\item $\nexists q\in\mathbb{Z}$ such that $\sqrt{\bar{\mu}}=(2q+1)\pi/(\overline{w}-\underline{w})$: To be able to satisfy $u(\overline{w})=u(\underline{w})=0$, it must be that $\alpha=0$. In this case, $\int_{w\in \mathcal{W}} u(w)^2\mathrm{d}w=0$, which contradicts the requirement that $\int_{w\in \mathcal{W}} \gamma(w)\mathrm{d}w=1$.
\item $\exists q\in\mathbb{Z}$ such that  $\sqrt{\bar{\mu}}=(2q+1)\pi/(\overline{w}-\underline{w})$: In this case, $\beta=(\overline{w}+\underline{w})/2$. To ensure that  $u(w)\neq 0$ for all $w\in\mathrm{int}(\mathcal{W})$, select $q=0$. Finally, to be able to satisfy the equality constraint $\int_{w\in \mathcal{W}} \gamma(w)\mathrm{d}w=1$, pick $\alpha=\sqrt{2/(\overline{w}-\underline{w})}$.
\end{itemize}
Finally, note that since the cost function and the constraint set are convex, all the extrema are minimizers.

\section{Proof of Corollary~\ref{cor:nonuniformp}} 
\label{proof:cor:nonuniformp}
It can be shown that~\eqref{eqn:tho:nonlinear_query} becomes
\begin{align*}
\trace\bigg(\begin{bmatrix}
1 & 1 \\
1 & 1
\end{bmatrix}
%\begin{bmatrix}
%\dfrac{\partial u(w,x)}{\partial w} \\[.5em]
%\dfrac{\partial u(w,x)}{\partial x}
%\end{bmatrix}
D^2 u(w,x)\bigg)
&+\frac{p'(x)}{p(x)} \bigg[\dfrac{\partial u(w,x)}{\partial w}+\dfrac{\partial u(w,x)}{\partial x}\bigg]\\
&+\mu(x)u(w,x)=0.
\end{align*}
Introducing the change of variable $v=x+w$ results in
\begin{align} \label{eqn:pdf:nonconstant_p}
\frac{\partial^2 u(v,x)}{\partial v^2} 
&+\frac{p'(x)}{p(x)} \dfrac{\partial u(v,x)}{\partial v}+\mu(x)u(v,x)=0.
\end{align}
Note that $v$ must belong to $[x+\underline{w},x+\overline{w}]$. To ensure that $u(v,x)=0$ on $x+\underline{w}$ and $x+\overline{w}$ (because $u(w,x)=0$ on $\partial \mathcal{W}$), $\mu(x)$ must be selected such that
\begin{align*}
\bigg[\frac{p'(x)}{p(x)}\bigg]^2-4\mu(x)<0, \forall x\in\supp(x).
\end{align*}
Under this condition, the solution of the partial differential equation in~\eqref{eqn:pdf:nonconstant_p} becomes
\begin{align*}
u(v,x)=&\alpha \exp\bigg(-\frac{p'(x)}{2p(x)}v\bigg)\\
&\times \cos\bigg(\sqrt{\bigg[\frac{p'(x)}{2p(x)}\bigg]^2-\mu(x)}(v-\beta) \bigg).
\end{align*}
Following the same line of reasoning as in the proof of Corollary~\ref{cor:problem2:averaging}, it can be inferred that $\beta=((x+\underline{w})+(x+\overline{w}))/2=x+(\underline{w}+\overline{w})/2$ and 
\begin{align*}
\sqrt{\bigg[\frac{p'(x)}{2p(x)}\bigg]^2-\mu(x)}=\frac{\pi}{\overline{w}-\underline{w}}.
\end{align*}
This concludes the proof.

\section{Proof of Corollary~\ref{cor:nonlinear=linear}}
\label{proof:cor:nonlinear=linear}
The proof follows from the selection of $u(w,x)=u(w)$ and $\mu(x)=\trace(F(x) F(x)^\top)\mu$ and following the same line of reasoning as in Corollary~\ref{cor:problem2:averaging}.\vspace{-.1in}

\section{Proof of Theorem~\ref{tho:nonlinear_query:unbounded}}
\label{proof:tho:4}\vspace{-.1in}
Similarly, because the cost function and the constraint set are convex, the stationarity condition is sufficient for optimality. Further, if multiple density functions satisfy the conditions, they all exhibit the same cost. In this case, the Lagrangian can be constructed as
\begin{align*}
\mathcal{L}
=&\int_{x\in\mathcal{X}}\int_{w\in \mathcal{W}}
\hspace{-.05in}p(x)\bigg(\hspace{-.07in}-\hspace{-.03in}\mu(x) \gamma(w|x)+\varrho w^\top w\gamma(w|x)\\
&\hspace{0.5in}+
\frac{1}{\gamma(w|x)}\bigg[F(x)^\top\frac{\partial \gamma(w|x)}{\partial w}+\frac{\partial \gamma(w|x)}{\partial x} \bigg]^\top\\
&\hspace{0.7in}\times \bigg[F(x)^\top\frac{\partial \gamma(w|x)}{\partial w}+\frac{\partial \gamma(w|x)}{\partial x} \bigg]\bigg) \mathrm{d}w\mathrm{d}x\\
&\hspace{0.5in}+\int_{x\in\mathcal{X}}\mu(x)p(x)\mathrm{d}x.
\end{align*}
The rest of the proof follows the same line of reasoning as in Theorem~\ref{tho:problem:general:2}.\vspace{-.1in}

\section{Proof of Corollary~\ref{cor:problem3:averaging}}
\label{proof:cor:problem3:averaging}\vspace{-.1in}
First, we present the solution of Problem~\ref{problem:not_meas:3}. 
In this proof, a solution of the form $u(w,x)=u(w)$ and $\mu(x)=\mu$ is sought for the partial differential equation in~\eqref{eqn:tho:nonlinear_query:unbounded}. Note that $u(w)=\alpha\exp(-w^\top \Sigma^{-1} w/4)$ satisfies the partial differential equation~\eqref{eqn:tho:nonlinear_query:unbounded} with $\mu=\trace(C^\top \Sigma^{-1} C)/2$ and $\Sigma=2(CC^\top)^{1/2}/\sqrt{\rho}$. Further to ensure that $\int_{w\in\mathcal{W}}u(w)^2=1$, select $\alpha=1/(\sqrt{(2\pi)^m\det(\Sigma)})$.  Now, we present the solution of Problem~\ref{problem:not_meas:3}. 
This follows from that for Problem~\ref{problem:not_meas:4} the duality gap is zero \cite{jeyakumar1990zero}. Therefore, the constraint on the variance can be added to the cost function using a Lagrange multiplier, which transforms the problem into that of Problem~\ref{problem:not_meas:3}. Therefore, following Corollary~\ref{cor:problem3:averaging}, the solution is equal to the density function of a zero-mean Gaussian random variable. Finally, for Gaussian random variables, the Fisher information is a decreasing function of the variance. Therefore, the Lagrange multiplier is set so that the inequality constraint on the variance becomes active. That means
$\trace(2(CC^\top)^{1/2}/\sqrt{\rho})=\vartheta.$ As a result, $\sqrt{\rho}=2\trace((CC^\top)^{1/2})/\vartheta$. \vspace{-.1in}

\section{Proof of Theorem~\ref{tho:dynamic_estimation}}
\label{proof:tho:dynamic_estimation}\vspace{-.1in}
Applying the result of Theorem~\ref{tho:nonlinear_query:unbounded} for linear query functions and probability density functions that are independent of the content of the server, the sufficient condition of optimality comes from the solution to the partial differential equation
$\trace(\Psi_T^\top D^2 u(w_T) \Psi_T)+(\mu-(\varrho/4) w_T^\top w_T) u(w_T)=0.$
Introduce the change of variable $w_T=\Psi_Tz$ for $z\in\mathbb{R}^n$ (recall that only one solution for the partial differential equation needs to be calculated).  It can be shown that $D^2\bar{u}(z)=D^2 u(\Psi_T z)=\Psi_T^\top D^2u(w_T)\big|_{w_T=\Psi_T z} \Psi_T$. Therefore,
$
\trace(D^2 u(\Psi_T z))+(\mu-(\varrho/4)  z^\top\Psi_T^\top \Psi_T z) u(\Psi_T z)=0.
$
Define $\bar{u}(z)=u(\Psi_T z)$. Thus,
$
\trace(D^2 \bar{u}(z))+(\mu-(\varrho/4)  z^\top\Psi_T^\top \Psi_T z) \bar{u}(z)=0.
$
The rest follows from that $\bar{u}(z)=\alpha\exp(-z^\top \Sigma^{-1} z/4)$ satisfies this partial differential equation. To do so, note that
$\trace(D^2 \bar{u}(z))=\bar{u}(z)\trace(\Sigma^{-1}zz^\top\Sigma^{-1})-\trace(\Sigma^{-1})/2.$ Therefore, $\mu=\trace(\Sigma^{-1})/2$, $\Sigma=2(\Psi_T^\top \Psi_T)^{-1/2}/\sqrt{\varrho}$, and $\alpha=1/\sqrt[4]{(2\pi)^m\det(\Sigma)}$.

\section{Proof of Proposition~\ref{prop:differential_privacy}}
\label{proof:prop:differential_privacy}
Note that
\begin{align*}
\frac{\exp\left(-|y-Cx|/b\right)}{\exp\left(-|y-Cx'|/b\right)}
&=\exp\left(|y-Cx'|/b-|y-Cx|/b\right)\\[-.4em]
&\leq \exp(|C(x-x')|/b)\\
&\leq \exp((\overline{x}-\underline{x})\max_i c_i/b)=\exp(\epsilon),
\end{align*}
where the first inequality follows from
\begin{align*}
|y-Cx'|
&=|(y-Cx')-(y-Cx)+(y-Cx)|\\
&\leq |(y-Cx')-(y-Cx)|+|y-Cx|\\
&\leq |C(x'-x)|+|y-Cx|.
\end{align*}
Therefore
$
\int_{y\in\overline{Y}}\exp\left(-|y-Cx|/b\right)\mathrm{d}y
\leq \exp(\epsilon)\linebreak\int_{y\in\overline{Y}}\exp\left(-|y-Cx'|/b\right)\mathrm{d}y.
$

\section{Proof of Proposition~\ref{prop:weak_differential_privacy}}
\label{proof:prop:weak_differential_privacy}
\vspace{-.1in}
Let $\mathrm{erfc}(x)$ denote the complementary error function defined as $\mathrm{erfc}(x):=\frac{2}{\sqrt{\pi}}\int_{x}^\infty \exp(-u^2/2)\mathrm{d}u.$
Define $\bar{K}:=\sqrt{2\log(1/(2\delta))}$. Evidently, $2\delta=\exp(-\bar{K}^2/2)$. From~\cite{chang2011chernoff}, it can be seen that $\delta\geq (1/2)\mathrm{erfc}(\bar{K}/\sqrt{2})$ and, therefore, $\sqrt{2}\mathrm{erfc}^{-1}(2\delta)\leq \bar{K}$. Note that
\begin{align*}
\vartheta
&\geq (\overline{x}-\underline{x})\max_{i}|c_i|\bigg(\frac{\sqrt{2\ln(1/(2\delta))}}{\epsilon}+\frac{1}{\sqrt{2\epsilon}} \bigg) \\[-.4em]
&=(\overline{x}-\underline{x})\max_{i}|c_i|\frac{1}{2\epsilon}\bigg(2\bar{K}+\sqrt{2\epsilon} \bigg)\\[-.4em]
&=(\overline{x}-\underline{x})\max_{i}|c_i|\frac{1}{2\epsilon}\bigg(\bar{K}+\sqrt{(\bar{K}+\sqrt{2\epsilon})^2} \bigg)\\[-.4em]
&\geq (\overline{x}-\underline{x})\max_{i}|c_i|\frac{1}{2\epsilon}\bigg(\bar{K}+\sqrt{\bar{K}^2+2\epsilon} \bigg).
\end{align*}
The above inequality in conjuction with~\cite{le2014differentially} shows that 
the Gaussian mechanism in Corollary~\ref{cor:problem3:averaging} is $(\epsilon,\delta)$-differentially private.

\end{document}